\documentclass[11pt]{article} 
\usepackage[utf8]{inputenc}
\usepackage{caption,graphicx,amsmath,amsthm,amssymb, authblk,xcolor,enumitem,subfigure,float,verbatim,hyperref,pgfplots,bibentry,mathtools,soul,isomath}
\usepackage{mathrsfs}
\usepackage{indentfirst}
\graphicspath{ {images/} }
\usepackage{multicol}
\usepackage{graphicx,amsmath,amssymb,amsthm,marvosym,tikz,fontenc}
\usetikzlibrary{backgrounds}
\usetikzlibrary{patterns,shapes.misc, positioning}
\usepackage{pgfplots}
\pgfplotsset{width=10cm,compat=1.9,tick scale binop=\times}
\usepackage{indentfirst}
\usepackage[a4paper,bindingoffset=0.2in,%
left=0.6in,right=0.8in,top=1.2in,bottom=1.2in,%
footskip=.25in]{geometry}
\usepackage{relsize}
\usepackage[english]{babel}
\allowdisplaybreaks
\theoremstyle{plain}
\newtheorem{theorem}{Theorem}[section]
\newtheorem{lemma}[theorem]{Lemma}

\newtheorem{corollary}[theorem]{Corollary}
\newtheorem{definition}[theorem]{Definition}
\newtheorem{remark}[theorem]{Remark}
\DeclareRobustCommand{\rchi}{{\mathpalette\irchi\relax}} %\rchi
\newcommand{\irchi}[2]{\raisebox{\depth}{$#1\chi$}}   %\rchi
\DeclareRobustCommand{\rgamma}{{\mathpalette\irgamma\relax}} %\rgmamma
\newcommand{\irgamma}[2]{\raisebox{\depth}{$#1\gamma$}}
\newlength{\defbaselineskip}
\newcommand{\setlinespacing}[1]%
{\setlength{\baselineskip}{#1 \defbaselineskip}}

\date{}
\begin{document}
\title{On the Adjacency and Seidel Spectra of Hypergraphs}
\author{Liya Jess Kurian$^1$\footnote{liyajess@gmail.com},  Chithra A. V$^1$\footnote{chithra@nitc.ac.in}
 \\ \small 
 1 Department of Mathematics, National Institute of Technology Calicut,\\\small
 Calicut-673 601, Kerala, India\\ \small	}
\maketitle
\thispagestyle{empty}
\begin{abstract}
  A hypergraph generalizes the concept of an ordinary graph. In an ordinary graph, edges connect pairs of vertices, whereas in a hypergraph, hyperedges can connect multiple vertices at a time. In this paper, we obtain a relationship between the characteristic polynomial of Seidel and adjacency matrices of hypergraph and also compute all the eigenvalues of some $k$-uniform hypergraphs. Moreover, we estimate the adjacency and Seidel spectra of the uniform double hyperstar and sunflower hypergraph. In addition to that, we determine the Seidel spectrum and main Seidel eigenvalues of hyperstar .  
\end{abstract}
\section{Introduction}
 Let $G^*=(V, E)$ be a hypergraph of order n with vertex set $V=\{v_1,v_2,\cdots,v_n\}$ and edge set $E=\{e_1,e_2,e_3,\cdots,e_m\}$, each hyperedge $e_i\in E$ is a subset of $V$ \cite{Bretto2013}. The rank of hypergraph $G^*$ is the maximum cardinality of its hyperedges, and co-rank is the minimum cardinality of its hyperedges. The order of hypergraph $G^*=(V, E)$ is the cardinality of $V$. The degree $d(v)$ of a vertex $v\in V$ is the number of hyperedges that contain $v$. A hypergraph $G^*$ is said to be $k$-uniform hypergraph \textnormal{\cite{Cooper2012,Kumar2017}} if the cardinality of each of its hyperedges is $k$ where $k\geq 2$. It is evident that an ordinary graph is a 2-uniform hypergraph. A hypergraph with $d(v_i)=r$ for all $v_i\in V$ is called an $r$-regular hypergraph. A hypergraph is said to be $(k,r)$-regular hypergraph if it is both  $k$-uniform  and $r$-regualr. The properties of $(k,r)$ regular hypergraph are studied in \cite{Kumar2017}. The adjacency matrix $A=(a_{ij})$ of $G^*$ \cite{Lin2017} is an $ n\times n $  matrix whose rows and columns are indexed by the vertices of $G^*$ and for all $ v_i,v_j \in V, $
 \begin{equation*}
a_{ij} =\left\{
   \begin{array}{ll}
      \mid \{e_k  \in E: \{v_i,v_j\} \subset  e_k\}\mid    &  \mbox{, } v_i \neq v_j, k=1,2,3,...,m   \\
       0  & \mbox{, } v_i=v_j
   \end{array}.
   \right.
\end{equation*} 
The adjacency spectrum of hypergraphs, in particular the generalized spectrum of power hypergraphs, are studied in \cite{Cardoso2020}. Let $G=(V',E')$ be an ordinary graph. Then, the power graph is formed by adding $(k-2)$ vertices to each edge of a graph $G$. Hyperstar can be considered as a power graph of a star graph. In \cite{Cardoso2022}, Cardoso  investigated hyperstars and their properties. The author also gave the adjacency spectrum of hyperstar.
\begin{theorem}\textnormal{\cite{Cardoso2022}}\label{Cardoso2022}
  The adjacency spectrum of hyperstar $S_{n}^k$ is $$ \sigma_{A}(S_{n}^k)= \begin{pmatrix}
-1 & k-2 & r_1 & r_2\\
(n-1)(k-2) & n-2 & 1 & 1
\end{pmatrix} $$ where $r_1$ and $r_2$ are the roots of the equation $\lambda^2-(k-2)\lambda-(n -1)(k -1) = 0$. 
\end{theorem}
Let $ J_{k,n}$ denote all one matrix of order $k\times n$ and $J_n$ and $I_n$ of order $n$ denote the all one and identity matrix, respectively. Then, the Seidel matrix $S$ of a hypergraph $G^*$ is defined as $S=J_n-I_n-2A$ \cite{Zakiyyah2022}. The matrices $A$ and $S$ of $G^*$ are real and symmetric. So, their eigenvalues are real.  For any square matrix $M$ we can find a scalar $\lambda$ such that $M\mathbf{x}=\lambda \mathbf{x} $ where $\mathbf{x}$ is the nonzero eigenvector corresponding to eigenvalue $\lambda$. Let $\lambda_1 \geq \lambda_2 \geq\cdots \geq \lambda_n$ and $\mu_1 \geq \mu_2 \geq \cdots \geq \mu_n$ are the eigenvalues of $A$ and $S$ respectively. The collection of all eigenvalues together with their multiplicities is known as the spectrum of $A\left(\text{or}~S\right)$ of $G^*$. Let $\lambda_{1},\lambda_{2},\lambda_{3},\cdots,\lambda_{d},$ be the distinct eigenvalues of an adjacency matrix $A$ of hypergraph $G^*$ with multiplicities $m_{1},m_{2},m_{3},...,m_{d}$. Then the adjacency spectrum of $G^*$ is denoted by,
$$ \sigma_{A}(G^*)=\begin{pmatrix}
               \lambda_{1} & \lambda_{2} & \lambda_{3} & \cdots & \lambda_{d}\\
               m_{1} & m_{2} & m_{3} & \cdots  & m_{d}
             \end{pmatrix}.$$
The Seidel energy $SE(G^*)$ of hypergraph $G^*$ is defined as the sum of the absolute values of the Seidel eigenvalues of $G^*$. In \textnormal{\cite{Cvetkovic1970}}, Cvetkovic proposed the idea of the main eigenvalue; an eigenvalue is said to be the main eigenvalue if it has an eigenvector in which the sum of the entries is not equal to zero, that is it has an eigenvector which is not orthogonal to $\mathbfit{j}$ where $\mathbfit{j}$ denotes a column vector whose all entries are equal to 1. Note that a Seidel eigenvalue of $S$ is said to be a main Seidel eigenvalue of $G^*$ if the eigenspace is not orthogonal to $\mathbfit{j}$. The tensor product of $n\times m$ matrix $M=(m_{ij})$ and $p \times q$ matrix $N$ is an $np \times mq$ matrix given by $(M\otimes N)_{(i,j)}= m_{ij}N$. Throughout $A$ and $S$ represents adjacency and Seidel matrix of the hypergraph $G^*$.\\
In this paper, we focus on the study of some classes of non-regular hypergraphs. In Section 2, we give basic definitions and results that will be used later. In Section 3, we determine the relationship between the characteristic polynomial of Seidel and the adjacency matrices of a hypergraph. Also, we obtain the seidel spectrum of $(k,r)-$regular hypergraph. In Section 4, the Seidel spectrum and main Seidel eigenvalues of hyperstar are calculated. Also, we estimate the Seidel energy of the hyperstar.  In section 5, we compute the adjacency spectrum and Seidel spectrum of uniform double hyperstar. In section 6, the adjacency and Seidel spectrum of the sunflower hypergraph are given.

	\section{Preliminaries}
 
This section gives basic definitions, terminologies, and facts used in the main results.
\begin{theorem}\textnormal{\cite{walkreg}}\label{numofwalks}
    Let $v_i$ and $v_j$ be two vertices of a hypergraph $G^{*}$. Then the number of walks of length $k$ from $v_i$ to $v_j$ of $G^*$ is the $(i,j)$\textsuperscript{th} entry of the matrix $A^k$.
\end{theorem}
\begin{definition}\textnormal{\cite{Cvetkovic1978}}\label{walkgenfn}
     The walk generating function of the number of walks of hypergraph $G^*$ is given by,
     $$ H_{G^*}(t)= \sum\limits_{l=0}^\infty N_l t^l $$ where $N_l$ denote the number of walks of length $l$ in $G^*$.
    
\end{definition}
\begin{theorem}\textnormal{\cite{Cvetkovic1970}}\label{thm2.4}
 Let $G$ be a multigraph of order $n$  and $A(G)$ be the adjacency matrix of $G$. If $\lambda_1, \lambda_2,\lambda_3,\cdots,\lambda_n$ be the eigenvalues of $A(G)$ corresponding to the mutually orthogonal normalized eigenvectors $x_1, x_2,$ $ x_3, \cdots,x_n$ and $X=(x_{ij}) $ be an orthogonal matrix of the eigenvectors of $A$. Then the total number of walks of length $l$ in $G$ is given by,
    \begin{equation*}
       N_l = \sum_{j=1}^n C_j\lambda_j^l,
   \end{equation*} where $\displaystyle  C_j={\Bigl( \sum_{i=1}^n x_{ij}\Bigr)^2}. $ 
\end{theorem}

\begin{definition}\textnormal{\cite{Cardoso2022}}
Let $S_{n}$ be a star with n vertices $\{ v_{0,0},$ $v_{1,1},$ $v_{2,1},$ $\cdots$ $ v_{n-1,1}\} $, then hyperstar $S_{n}^k=(V,E)$ is obtained from the star by adding $k-2$ new vertices to each hyperedge in such a way that $V=\{ v_{0,0},v_{1,1},v_{1,2},$ $\cdots$ $,v_{1,k-1},v_{2,1},$ $v_{2,2},\cdots ,$ $ v_{2,k-1},\cdots ,$ $ v_{n-1,k-1}\}$ and $n-1$ hyperedges $E=\{\{v_{0,0},v_{1,1},v_{1,2},\cdots ,v_{1,k-1}\},$  $\{v_{0,0},v_{2,1},$ $v_{2,2},\cdots,$ $v_{2,k-1}\},$ $\cdots, $ $\{v_{0,0},v_{n-1,1},v_{n-1,2}, \cdots,v_{n-1,k-1}\}\}$. 
\end{definition}
\begin{figure}[htbp]
\centering
\tiny
\begin{tikzpicture}[scale=0.5]
\draw[ draw= black] (-3.5,-3.5) ellipse (1cm and 3cm);
\draw[ draw= black,rotate=270] (6,-1) ellipse (1cm and 3cm);
\draw[ draw= black,rotate=270] (6,-6) ellipse (1cm and 3cm);
\filldraw[fill=black](-3.5,-6)circle(0.1cm);  
\node at(-3.5,-7){$v_{00}$};
\filldraw[fill=black](-5.5,-6)circle(0.1cm);
\filldraw[fill=black](-7.5,-6)circle(0.1cm);  
\node at(-5.5,-6.5){$v_{11}$};
\node at(-7.5,-6.5){$v_{12}$};
\filldraw[fill=black](-1.5,-6)circle(0.1cm);
\filldraw[fill=black](0.5,-6)circle(0.1cm);  
\node at(-1.5,-6.5){$v_{21}$};
\node at(0.5,-6.5){$v_{22}$};
\filldraw[fill=black](-3.5,-4)circle(0.1cm);
\filldraw[fill=black](-3.5,-2)circle(0.1cm);  
\node at(-3.5,-3.5){$v_{31}$};
\node at(-3.5,-1.5){$v_{32}$};
\end{tikzpicture}
\caption{$S_4^3$-Hyperstar}
\label{fig:1}
\end{figure}
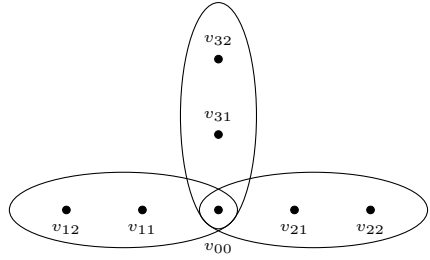

\begin{definition}\textnormal{\cite{Berge1973}}
 The complete $r$-uniform hypergraph $K_n^r $ is a hypergraph with $n$ vertices such that all possible subsets with $r$ vertices form hyperedges. 
\end{definition}

  \begin{lemma}\label{detblock}\textnormal{\cite{Das2018}}
             Let $B, C, W$, and $X$ be matrices with $B$ invertible. Let 
             $$ M= \begin{pmatrix}
                 B & C\\
                 D & X
                  \end{pmatrix}$$
                  Then $det(M)=det(B)det(X-DB^{-1}C)$ and if $X$ is invertible, then $det(M) = det(X)det(B-CX^{-1}D).$
          \end{lemma}
\begin{lemma}\textnormal{\cite{vrabel2016}}\label{vrabel2016}
Let $\mathbf{M}\in\mathbb{R}^{n\times n}$ be an invertible  matrix, and $U$ and $W$ are $n \times 1$ matrices. Then $$\text{det}(\mathbf{M}+UW^{T})=\text{det}(\mathbf{M})+W^{T}\text{adj}(\mathbf{M})U$$
where $adj(\mathbf{M})$ denotes the adjoint of $\mathbf{M}.$
\end{lemma}
\begin{definition}\textnormal{\cite{supertrees}}
    Let $S_{n_1,n_2}$ be a double star of order $n_1+n_2$, which is obtained by adding an edge connecting the central vertices of the star $S_{n_1}$ and $S_{n_2}$. Then the $k$-th power of $S_{n_1,n_2}$ is called uniform double hyperstar $S_{n_1,n_2}^k.$ 
\end{definition}

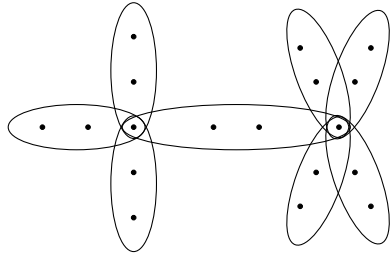
\begin{figure}[htbp]
\centering
\tiny
\begin{tikzpicture}[scale=0.3]
\draw[ draw= black,rotate=90] (0,0) ellipse (1cm and 5cm);
\draw[ draw= black] (-4.5,2.5) ellipse (1cm and 3cm);
\draw[ draw= black] (-4.5,-2.5) ellipse (1cm and 3cm);
\draw[ draw= black,rotate=90] (0,7) ellipse (1cm and 3cm);
\draw[ draw= black,rotate=200] (-4.2,-1) ellipse (1cm and 3cm);
\draw[ draw= black,rotate=160] (-4.2,1) ellipse (1cm and 3cm);

\draw[ draw= black,rotate=20] (4.2,-4) ellipse (1cm and 3cm);
\draw[ draw= black,rotate=340] (4.2,4) ellipse (1cm and 3cm);
\filldraw[fill=black](1,0)circle(0.1cm);
\filldraw[fill=black](-1,0)circle(0.1cm);
\filldraw[fill=black](-4.5,0)circle(0.1cm);
\filldraw[fill=black](4.5,0)circle(0.1cm);
\filldraw[fill=black](-4.5,4)circle(0.1cm);
\filldraw[fill=black](-4.5,2)circle(0.1cm);
\filldraw[fill=black](-4.5,-4)circle(0.1cm);
\filldraw[fill=black](-4.5,-2)circle(0.1cm);
\filldraw[fill=black](-6.5,0)circle(0.1cm);
\filldraw[fill=black](-8.5,0)circle(0.1cm);
\filldraw[fill=black](5.2,2)circle(0.1cm);
\filldraw[fill=black](5.9,3.5)circle(0.1cm);
\filldraw[fill=black](5.2,-2)circle(0.1cm);
\filldraw[fill=black](5.9,-3.5)circle(0.1cm);
\filldraw[fill=black](3.5,2)circle(0.1cm);
\filldraw[fill=black](2.8,3.5)circle(0.1cm);
\filldraw[fill=black](3.5,-2)circle(0.1cm);
\filldraw[fill=black](2.8,-3.5)circle(0.1cm);
\end{tikzpicture}
\caption{Uniform double hyperstar-$S_{4,5}^3$}
\label{fig:my3}
\end{figure}

\begin{theorem}\textnormal{\cite{haemers1995interlacing}}\label{interlacingthm}
Let $N\in \mathbb{R}^{n\times n}$ be a symmetric matrix with eigenvalues $\lambda_1\geq\lambda_2\geq\lambda_3\geq\cdots\geq\lambda_n$. If $\mu_1\geq\mu_2\geq\mu_3\geq\cdots\geq\mu_m$ are the eigenvalues of the  principal submatrix $M\in \mathbb{R}^{m\times m}$, then 
    % Let $N$ be a symmetric square matrix of order $n$ with eigenvalues $\lambda_1\geq\lambda_2\geq\lambda_3\geq\cdots\geq\lambda_n$. If $M$ is a principal submatrix of order $m$ with eigenvalues $\mu_1\geq\mu_2\geq\mu_3\geq\cdots\geq\mu_m$, then the eigenvalues of $M$ interlace eigenvalues of $N$,
    $$\lambda_i\geq\mu_i\geq\lambda_{n-m+i}~\text{for}~i=1,2,\cdots,m .$$
\end{theorem}

  \begin{lemma}\textnormal{\cite{jahfar2020}}\label{detij}
              For any two real numbers $r$ and $s$,
              $$(rI_n-sJ_n)^{-1}=\frac{1}{r}I_n+\frac{s}{r(r-ns)}J_n.$$
          \end{lemma}
    \begin{lemma}\textnormal{\cite{jahfar2020}}\label{coronal}
              Let $M$ be an $n\times n $ matrix.Then
              $$det(\beta I_n-M-\rgamma J_n)=(1-\rgamma \rchi_M(\beta))det(\beta I_n-M).$$
              For an $n\times n$ real matrix $M$ with row sum equal to $r$, $\rchi_M(\beta)=\frac{n}{\beta-r}$
          \end{lemma}
% \begin{lemma}\textnormal{\cite{Das2018}}
%     Let $M\in \mathbb{R}^{n\times n}$ with each row sum equal to a constant $c$, then $\rchi_M(x)=\frac{n}{x-c}$
% \end{lemma}

\begin{definition}\textnormal{ \cite{Hu2013}}
          Let $S^k=(V,E)$ be a $k$-uniform hypergraph of order $k(k-1)+1$. If label the vertex set $V$ as $V=\{v_{0,0},v_{1,1},v_{1,2},\cdots , v_{1,k}, \cdots, v_{k-1,1},v_{k-1,2}, v_{k-1,3}, \cdots, v_{k-1,k}\}$  such that set of hyperedges being $E=\{\{v_{1,1},v_{1,2},v_{1,3} \cdots , v_{1,k}\},\cdots ,\{v_{k-1,1}, v_{k-1,2}, v_{k-1,3},\cdots ,v_{k-1,k}\},\{v_{0,0}, v_{1,1},v_{2,1}, \cdots,$ $ v_{k-1,1}\}\}$, then $S^k$ is a sunflower hypergraph.
      \end{definition}   
         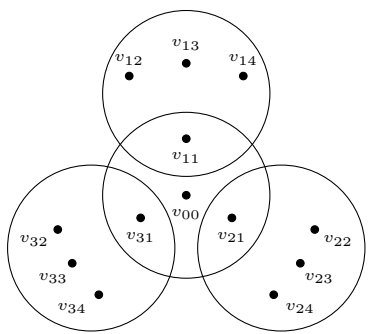
\begin{figure}[H]
\centering
\tiny
\begin{tikzpicture}[scale=0.5]
\draw(0,2.7)circle(2.2cm);
\draw(0,0)circle(2.2cm);
\draw(2.5,-1.4)circle(2.2cm);
\draw(-2.5,-1.4)circle(2.2cm);
% \draw(0,0)circle(3.5cm);
\filldraw[fill=black](0,3.5)circle(0.1cm);
\filldraw[fill=black](-1.5,3.162)circle(0.1cm);
\filldraw[fill=black](1.5,3.162)circle(0.1cm);
\filldraw[fill=black](-2.3,-2.638)circle(0.1cm);
\filldraw[fill=black](2.3,-2.638)circle(0.1cm);
\filldraw[fill=black](-3,-1.803)circle(0.1cm);
\filldraw[fill=black](3,-1.803)circle(0.1cm);
\filldraw[fill=black](-3.38,-0.909)circle(0.1cm);
\filldraw[fill=black](3.38,-0.909)circle(0.1cm);
\filldraw[fill=black](0,0)circle(0.1cm);
\filldraw[fill=black](0,1.5)circle(0.1cm);
\filldraw[fill=black](-1.2,-0.6)circle(0.1cm);
\filldraw[fill=black](1.2,-0.6)circle(0.1cm);
\node at (0,-0.5) {$v_{00}$};
\node at (0,1) {$v_{11}$};
\node at (1.2,-1.1) {$v_{21}$};
\node at (-1.2,-1.1) {$v_{31}$};
\node at (3,-3) {$v_{24}$};
\node at (3.5,-2.2) {$v_{23}$};
\node at (4,-1.2) {$v_{22}$};

\node at (-3,-3) {$v_{34}$};
\node at (-3.5,-2.2) {$v_{33}$};
\node at (-4,-1.2) {$v_{32}$};
\node at (1.5,3.6) {$v_{14}$};
\node at (-1.5,3.6) {$v_{12}$};
\node at (0,4) {$v_{13}$};
\end{tikzpicture}
\caption{sunflower hypergraph: $S^4$}
\label{fig:my2}
\end{figure}

          \begin{theorem}\textnormal{\cite{Hagos2002}}\label{Hagos2002}
              Let $G$ be a graph of order $n$. Then the rank of the matrix $\begin{bmatrix}
                  \mathbfit{j} & A\mathbfit{j} &\cdots& A^{n-1}\mathbfit{j}
              \end{bmatrix}$  is equal to the number of its main eigenvalues. 
          \end{theorem}
          Let $M$ be a real matrix of order $n$  such that rows and columns of $M$ are indexed by elements of $X=\{1,2,3,\cdots,n\}$. Consider a partition $P=\{X_1,X_2,\cdots,X_m\}$ of $X$. Then the partition of $M$ according to $P$ is $\begin{bmatrix}
              M_{11} & M_{12} & \cdots & M_{1m}\\
              M_{21} & M_{22} & \cdots & M_{2m}\\
              \vdots  & \vdots & \ddots & \vdots\\
              M_{m1} & M_{m2} & \cdots & M_{mm}
          \end{bmatrix},$ where each $M_{ij}$ is a submatrix  of $M$ such that rows and columns of $M_{ij}$ are indexed by elements of $X_i$ and $X_j$ respectively. If $q_{ij}$ denotes the average row sum of $M_{ij}$, then the matrix $Q= (q_{ij})$ is called a quotient matrix of $M$.  If the row sum of each block $M_{ij}$ is a constant, then the partition $P$ is called equitable.   
          \begin{theorem}\label{equipartition}\textnormal{\cite{Fouzul2018}}
               Let $Q$ be a  quotient matrix of any square matrix $M$ corresponding to an equitable partition.  Then, the spectrum of $M$ contains the spectrum of $Q$.
          \end{theorem}
         
 \section{Characteristic Polynomial of a Hypergraph}
	Using the characteristic polynomial of the adjacency matrix of hypergraph, we can find the spectrum of the Seidel matrix of  $G^*$.\\
 The following theorem gives the relation between the characteristic polynomial of the adjacency and the Seidel matrix of the hypergraph $G^*$.
	\begin{theorem}\label{cpas}
         Let $P_S(\lambda)$ be the characteristic polynomial of the Seidel matrix of $G^*$ and $P_A(\lambda)$ be the characteristic polynomial of the adjacency matrix of $G^*$. Then, 
		\begin{equation*}
			P_S(\lambda)=(-2)^nP_A\Bigl(-\frac{\lambda+1}{2}\Bigr) \left(\frac{-1}{\lambda+1}H_{G^*}\Bigl(\frac{-2}{\lambda+1}\Bigr)+1\right),
		\end{equation*}
   where $H_{G^*}$ is the walk generating function of number of walks in $G^*$. 
	\end{theorem}
	\begin{proof}
		Consider an invertible square matrix $\mathbf{M}$, $Sum(\mathbf{M})$ denote the sum of all of the entries of $\mathbf{M}$. 
  % \hl{$adj(\mathbf{M})$ denotes the adjoint matrix for $\mathbf{M}$}.
  % We know that $adj(\mathbf{M})=\mathbf{M}^{-1}det(\mathbf{M})$. 
  By Lemma \ref{vrabel2016},
 % Now consider the column matrices $U$ and $W$, then by \hl{ref from preliminary-}matrix determinant Lemma %\cite{vrabel2016},
		$$\text{det}(\mathbf{M}+UW^{T})=\text{det}(\mathbf{M})+W^{T}\text{adj}(\mathbf{M})U.$$
		Take $U_{n\times1}=[\frac{1}{2} \;\;  \frac{1}{2}  \;\;  \frac{1}{2}\; \cdots \; \frac{1}{2}]^{T}$ and $W_{n\times 1}=[w  \;\; w \;\;  w\; \cdots\; w]^{T}$, where $w\in \mathbb{R}$ then
                $$ \; UW^{T}=\frac{w}{2}\;J \; \; \; \;\text{and} \; \; \; \; W^{T}\text{adj}(\mathbf{M})U=\frac{w}{2} Sum(adj(\mathbf{M})),$$ So,
		\begin{equation}\label{cp2}
			\text{det}(\mathbf{M}+\frac{w}{2}J)=\text{det}(\mathbf{M})+\frac{w}{2} Sum(adj(\mathbf{M})).
		\end{equation}
		From Theorem \ref{numofwalks}, the total number of walks of length $l$, $N_l= \sum\limits_{i,j=1}^{n} a_{ij}^l$, where $a_{ij}^l$ is the $ij$-th entry of $A^l$. Therefore,$N_l=Sum(A^l).$ Let $ H_{G^*}(t)= \sum\limits_{l=0}^\infty N_l t^l $ be the generating function of the number of walks of length $l$ of $G^*$. Then,
		\begin{equation*}\label{equ3.14}
		    H_{G^*}(t)= \sum\limits_{l=0}^\infty N_l t^l =\sum_{l=0}^\infty Sum(A^l) t^l.
		\end{equation*}
 We know that, $	\displaystyle\sum_{l=0}^\infty A^l t^l =(I-tA)^{-1}\; \text{when} ~\displaystyle\| tA \| < 1.$
Then, $\displaystyle\sum_{l=0}^\infty A^l t^l =\text({det}(I-tA))^{-1}adj(I-tA).$	\\
% We have,	
% $$\sum_{l=0}^\infty Sum(A^l) t^l=\left(\text{det}(I-tA)\right)^{-1}Sum(adj(I-tA)).$$
 Therefore,
  $$ H_{G^*}(t)=\left(\text{det}(I-tA)\right)^{-1}Sum(adj(I-tA)).$$
	From (\ref{cp2}) we get,
		$Sum(adj(I-tA))=\displaystyle\frac{2}{t} [\text{det}(I-tA+\frac{t}{2}J)-\text{det}(I-tA)].$\\
 Thus,
		\begin{equation*}\label{cp4}
			Sum(adj(I-tA))=\frac{2}{t} [\text{det}((1+\frac{t}{2})I+\frac{t}{2}S)-\text{det}(I-tA)].
		\end{equation*}
Hence $H_{G^*}(t)$ ,
		\begin{align*}
			H_{G^*}(t) =\frac{2}{t}\left(\frac{\text{det}((1+\frac{t}{2})I+\frac{t}{2}S)-\text{det}(I-tA)}{\text{det}(I-tA)}\right)
			&=\frac{2}{t}\left(\frac{\text{det}(\frac{1}{2}(\frac{2+t}{t})I+S)}{\text{det}(\frac{1}{t}I-A)}-1\right)\\
			&=\frac{2}{t}\left(\left(\frac{-1}{2}\right)^n\frac{\text{det}(-(\frac{2+t}{t})I-S)}{\text{det}(\frac{1}{t}I-A)}-1\right).
		\end{align*}
Therefore,
		\begin{equation*}
			H_{G^*}(t) =\frac{2}{t}\left(\left(-\frac{1}{2}\right)^n\frac{P_S\left( -\frac{t+2}{t}\right)}{P_A(\frac{1}{t})}-1\right).
		\end{equation*}
            Then 
		\begin{equation}\label{cp5}
			 P_A(\lambda)=\frac{\left( -\frac{1}{2}\right)^nP_S(-1-2\lambda)}{\frac{1}{2\lambda}H_{G^*}\left(\frac{1}{\lambda}\right)+1},~\text{when}~ t=\displaystyle{\frac{1}{\lambda}}
		\end{equation}
  and replacing $t$ by $\frac{-2}{1+\lambda}$
  \begin{equation*}
			 P_S(\lambda)=(-2)^nP_A\Bigl(-\frac{\lambda+1}{2}\Bigr) \left(\frac{-1}{\lambda+1}H_{G^*}\Bigl(\frac{-2}{\lambda+1}\Bigr)+1\right).
		\end{equation*}
  
	\end{proof}
	\begin{lemma}\label{lemma3.2}
 Let $G^*$ be a hypergraph of order $n$ and $X=(x_{ij}) $ be a matrix of mutually orthogonal normalized eigenvectors of $A$ corresponding to the eigenvalues $\lambda_1, \lambda_2,\lambda_3,\cdots,\lambda_n$. Then the total number of walks of length $l$ in $G^*$ is given by,
    \begin{equation*}
       N_l = \sum_{j=1}^n C_j\lambda_j^l,
   \end{equation*} where $\displaystyle  C_j={\Bigl( \sum_{i=1}^n x_{ij}\Bigr)^2}. $ 
	\end{lemma}
 \begin{proof}
 The proof follows from Theorem \ref{thm2.4}.
     
 \end{proof}
	\begin{theorem}\label{seideleig}
	    If the adjacency spectrum of hypergraph $G^*$ contains an eigenvalue $\lambda_0 $ with multiplicity $m_p>1$, then the  Seidel spectrum of $G^*$ has an eigenvalue $-2\lambda_0-1$ with multiplicity $m_q$, where $m_q\geq m_{p-1}$.
	\end{theorem}
	\begin{proof}
By Definition \ref{walkgenfn} and Lemma \ref{lemma3.2},
		
			$$H_{G^*}(t) =\sum_{j=1}^n C_j \frac{1}{1-t\lambda_j}.$$
Now we define the function $\Phi$ as, 
                     \begin{equation*}\label{eqn12.1}
                         \Phi(\lambda)=\frac{(-\frac{1}{2})^n P_S(-1-2\lambda)}{P_A(\lambda)}.
                     \end{equation*}
From  (\ref{cp5})  we obtain,
		\begin{equation*}
			\Phi(\lambda) =\frac{1}{2\lambda}H_{G^*}(\frac{1}{\lambda})+1=\frac{1}{2}\sum_{j=1}^n C_j\frac{\lambda}{\lambda-\lambda_j}+1.
		\end{equation*}
By expanding the right hand side of the above equation, we get a rational polynomial $\displaystyle \frac{P_1(\lambda)}{P_2(\lambda)} $. Since
                $$ \frac{(-\frac{1}{2})^nP_S(-1-2\lambda)}{P_A(\lambda)}=\frac{1}{2}\sum_{j=1}^n C_j\frac{\lambda}{\lambda-\lambda_j}+1,$$
		%where $P_1$ and $P_2$ are polynomials in $\lambda$.
  it is clear that the roots of $P_2(\lambda)$ are all of multiplicity 1. So if $\lambda_0$ is an eigenvalue of $A$ with multiplicity $m_p\:(m_p\geq 2)$, $P_S(-1-2\lambda)$ contain a factor $(\lambda-\lambda_0)^{m_q}$ where $m_q\geq m_{p-1}$. Therefore, $ P_S(\lambda)$ contains the factor $(\lambda+2\lambda_0+1)^{m_q}$.
		Hence, Seidel spectrum of $G^*$ contains an eigenvalue $-2\lambda_0 -1 $ with multiplicity $m_q$.
	\end{proof}
	\subsection{Characteristic polynomial of \texorpdfstring{$(k,r)$}{e}-regular hypergraph}
	The $(k,r)$ regular hypergraph was investigated in \cite{Kumar2017,Li1996} . A $(k,r)$-regular hypergraph is a $k$-uniform $r$-regular hypergraph. In \cite{Li1996}, Li W and Solé P derive $ r(k-1)$ as an eigenvalue of the adjacency matrix of the $(k,r)$ regular hypergraph. 
	\begin{lemma}\label{Nl}
		Let $G^*$ be a $(k,r)$-regular hypergraph with vertices $v_i\;(1\leq i\leq n)$ and hyperedges $e_j\;(1\leq j\leq m)$. Then the number of walks of length $l$ is given by, $$N_l=nr^l(k-1)^l.$$
	\end{lemma}
	
	\begin{proof}
 For a $(k,r)$-regular hypergraph $G^*$ on $n$ vertices, let $N_l$ be the number of walks of length $l$. The proof follows from induction on length $l$. When $l=1$, pick any random vertex $v_1\in G^*$, and assume that it is contained in hyperedge $e_1$, which contains $k-1$ additional vertices. By this argument, there are $k-1$ walks of length one with origin $v_1$. Since $G^*$ is $r$-regular, there exist $r(k-1)$ walks of length $1$ starting from $v_1$. Hence,
		$$ N_1=nr(k-1).$$
  
		% For  $l=2$, consider any arbitrary walk of length one, say $v_1e_1v_2$. There are $r(k-1)$ walks of length one starting from $v_1$. Since  $v_1e_1v_2$ is an arbitrary walk starting from $v_1$, the possible number of such walks is $nr(k-1)$. Therefore,
		% $$ N_2=nr^2(k-1)^2.$$ 
     Assume that the result holds for $l=p$ then $N_p= nr^p(k-1)^p.$ Now, we prove for $l=p+1$. For that we choose a walk of length $p$ from the $nr^p(k-1)^p$ walks, say $v_1e_1v_2e_2$ $v_3\cdots e_{p-1}v_{p}$. Since $v_{p}$ is adjacent to $r(k-1)$ vertices, we can obtain $r(k-1)$ walks. The walk $v_1e_1v_2e_2$ $v_3\cdots e_{p-1}v_{p}$ is arbitary, so in total we can have $(nr^p(k-1)^p)(r(k-l))=nr^{p+1}(k-1)^{p+1}$ walks of length $l+1$. Hence the theorem.
	\end{proof}
	\begin{lemma}\label{HF}
		The generating function of the number of walks of a $(k,r)$-regular hypergraph $G^*$ on n vertices is given by,
		$$H_{G^*}(t)=\frac{n}{1-r(k-1)t}\;\;\; \text{if}\;|t|\leq\frac{1}{r(k-1)}.$$
	\end{lemma}
	\begin{proof}
 The proof follows from Lemma \ref{Nl}
 \begin{comment}
     The generating function of the number of walks of $G^*$ is defined as,
		$$H_{G^*}(t)=\sum_{l=0}^\infty N_lt^l.$$
		By Lemma \ref{Nl} we obtain,
		\begin{align*}
			H_{G^*}(t) =\sum_{l=0}^\infty nr^l(k-1)^l t^l=\frac{n}{1-r(k-1)t}.
		\end{align*}
 \end{comment}
		
	\end{proof}
	\begin{theorem}\label{Seidel}
		Let $G^*$ be a $(k,r)$-regular hypergraph with $n$ vertices and it's adjacency spectrum is $\lambda_1=r(k-1),\lambda_2,\lambda_3,\cdots,\lambda_n $. Then the Seidel spectrum of  $G^*$ is $ n-1-2\lambda_1, -1-2\lambda_2,-1-2\lambda_3,\cdots,-1-2\lambda_n$.
	\end{theorem}
	\begin{proof}
		From Theorem \ref{cpas} and Lemma \ref{HF} we obtain,
		$$\frac{n}{1-r(k-1)t}=\frac{2}{t}\left(\left(-\frac{1}{2}\right)^n\frac{P_S\left( -\frac{t+2}{t}\right)}{P_A(\frac{1}{t})}-1\right).$$
		Putting $-\displaystyle{\left(\frac{t+2}{t}\right)=\lambda }$ we have,
		\begin{align*}
			-\lambda-1\left( \left(-\frac{1}{2}\right)^n \frac{P_S(\lambda)}{P_A(\frac{-\lambda-1}{2})}-1\right)=\frac{n}{1-r(k-1)\frac{-2}{(\lambda+1)}}
			=\frac{n(\lambda+1)}{(\lambda+1)+2r(k-1)}
		\end{align*}
		On simplification we get,
		\begin{equation*}
			P_S(\lambda)=(-2)^n\left[ \frac{-n+\lambda+1+2r(k-1)}{\lambda+1+2r(k-1)}\right]P_A\left(\frac{-\lambda-1}{2}\right).
		\end{equation*}
		Since $r(k-1)$ is an eigenvalue of a $(k,r)$-regular hypergraph, $P_A(\lambda)$ contains the factor $(\lambda-r(k-1))$. Thus,
		$$ P_A\Bigl(\frac{-\lambda-1}{2}\Bigr)=(\lambda+1+2r(k-1))Q_A\Bigl(\frac{-\lambda-1}{2}\Bigr),$$
		where $Q_A$ is another polynomial of degree one less than $P_A$. Therefore,
		$$ P_S(\lambda)=(-2)^n (\lambda -n+1+2r(k-1))Q_A\Bigl(\frac{-\lambda-1}{2}\Bigr).$$
		The Seidel eigenvalue corresponding to the eigenvalue $r(k-1)$ of the adjacency spectrum is $n-1-2r(k-1)$. Since $G^*$ is a $(k,r)$ regular hypergraph with n vertices,
		$$P_S(\lambda)=(-2)^n \frac{\lambda-n+1+2r(k-1)}{\lambda+1+2r(k-1)}P_A\Bigl(\frac{-\lambda-1}{2}\Bigr).$$
	\end{proof}
	An $r$-uniform complete hypergraph of order $n$, as referred by Berge in \cite{Berge1989}, is a hypergraph consisting of all the $r$-subsets of the vertex set $V$. Zakiyyah \cite{Zakiyyah2022} deals with the spectrum of $r$-uniform complete hypergraph.  If the adjacency spectrum is known, we use Theorem \ref{Seidel} to establish the Seidel spectrum of the hypergraph.\\
	For example, the adjacency spectrum of $K_n^r$ is,

  \[ \sigma_{A}(K_n^r)= \begin{pmatrix}
(n-1) \left( \begin{smallmatrix} n-2\\
			r-2
		\end{smallmatrix}\right) & -\left( \begin{smallmatrix} n-2\\
			r-2
		\end{smallmatrix}\right) \\
1 & ~~n-1 
\end{pmatrix}. \]

\noindent Using Theorem \ref{Seidel}, the Seidel spectrum of $(K_{n}^{r})$ is
		 \[ \sigma_{S}(K_n^r)=\begin{pmatrix}
(n-1)(1-2 \left( \begin{smallmatrix} n-2\\
			r-2
		\end{smallmatrix}\right)) & 2\left( \begin{smallmatrix} n-2\\
			r-2
		\end{smallmatrix}\right)-1 \\
1 & n-1 
\end{pmatrix}. \]

 \section{Seidel spectrum of hyperstars}
 
This section extends the study to the Seidel spectrum of hyperstars. Hyperstars are $k$-uniform hypergraphs with $(n-1)(k-1)+1$ vertices and $n-1$ hyperedges. \\

Let $G^*=(V, E)$ be the hypergraph with n vertices and $\mathbf{x}=(x_{v_i})$ be an $n$-dimensional vector. Define $$x(\alpha)=x(v_1,v_2,v_3,\cdots, v_t)=x_{v_1}+x_{v_2}+x_{v_3}+\cdots+x_{v_t},$$  where $\alpha$ is the non-empty collection of subset of $V$. For simplicity, we write $x_{v_i}$ as $x_{i}$. Let $E_{[v]}$ be the set of all hyperedges containing vertex $v$. The entries corresponding to the vertex $v$ of the adjacency matrix $A$ of $G^*$ is given by $(A)_v$.
Then,
\begin{equation*}
    (A\mathbf{x})_v=\sum_{e\in E_{[v]}}x(e-\{v\}), \forall v\in V.
\end{equation*}
Then for a Seidel matrix $S$ of $G^*$,
\begin{equation}\label{1}
    (S\mathbf{x})_v=x(V-\{v\})-2\sum_{e\in E_{[v]}}x(e-\{v\}), \forall v\in V.
\end{equation}

\noindent For example,
let $G^*$ be a hypergraph with vertex set $V=\{v_1,v_2,v_3,v_4,v_5\}$ and hyperedge set $E=\{\{v_1,v_2,v_3\},$ $\{v_2,v_3,$ $v_4,v_5\},$ $\{v_1,v_2,v_4\}\}$ and $\mathbf{x}=\begin{bmatrix}
    x_1 & x_2 & x_3 & x_4 & x_5
\end{bmatrix}^{T}$. Then,
$$S\mathbf{x}=\begin{bmatrix}
    0 & -3  & -1 & -1 & 1\\
    -3 & 0  & -3 & -3 & -1\\
    -1 & -3 & 0 & -1 & -1\\
    -1 & -3 & -1 & 0 & -1\\
    1 & -1 & -1 & -1 & 0
\end{bmatrix}\begin{bmatrix}
    x_1\\
    x_2\\
    x_3\\
    x_4\\
    x_5
\end{bmatrix}=\begin{bmatrix}
    -3x_2-x_3-x_4+x_5\\
    -3x_1+-3x_3-3x_4-x_5\\
    -x_1-3x_2-x_4-x_5\\
    -x_1-3x_2-x_3-x_5\\
    x_1-x_2-x_3-x_4
\end{bmatrix}.$$

\noindent Therefore, $(S\mathbf{x})_{v_1}=-3x_2-x_3-x_4+x_5$. Since $E_{[v_1]}=\{\{v_1,v_2,v_3\},$ $\{v_1,v_2,v_4\}\}$, we get
$$\sum_{e\in E_{[v_1]}}x(e-\{v_1\})=x(v_2,v_3)+x(v_2,v_4)=2x_2+x_3+x_4.$$
Also,
$$x(V-\{v_1\})=x(v_1,v_2,v_4,v_5)=x_2+x_3+x_4+x_5.$$
Therefore from (\ref{1}) we get,
$$(S\mathbf{x})_{v_1}=x_2+x_3+x_4+x_5-2(2x_2+x_3+x_4)=-3x_2-x_3-x_4+x_5.$$

	\begin{lemma}\label{xu_equals_xv}
 	Let $G^*$ be a $k$-uniform hypergraph and $u$, $v \in V(G^*)$, belonging to the exact same hyperedges. If $(\lambda,\mathbf{x})$ is an eigenpair of  $S$ with  $\lambda+1 \neq 2d(u) $ then $x_u=x_v$, where $x_u$ and $x_v$ are entries of $\mathbf{x}$ corresponding to the vertices $u$ and $v$ respectively.
	\end{lemma}
	\begin{proof}
		Let $(\lambda,\mathbf{x})$ be an eigenpair of $S(G^*)$ and $\mathbf{x}=(x_1, x_2, x_3,\cdots,x_n)$  be the corresponding eigenvector.\\
		By the definition of $S$,
		\begin{equation*}
			((S+I)\mathbf{x})_u=((J-2A)\mathbf{x})_u.
		\end{equation*}
		Now,
		\begin{equation*}
			\begin{split}
				((S+I)\mathbf{x})_u-2d(u)x_u & = ((J-2A)\mathbf{x})_u-2d(u)x_u\\
				& = x_1+x_2+x_3+...+x_n-2(\sum_{e\in E_{[u]}} x(e-\{u\})+d(u)x_u)\\
				& = x_1+x_2+x_3+...+x_n-2\sum_{e\in E_{[u]}} x(e)\\
				& = x_1+x_2+x_3+...+x_n-2\sum_{e\in E_{[v]}} x(e)\\
				& = ((S+I)\mathbf{x})_v-2d(v)x_v.
			\end{split}
		\end{equation*}
		Since $(\lambda,\mathbf{x})$ is an eigenpair of $S$, $( \lambda +1,\mathbf{x})$ is an eigenpair of $S+I$.\\
		Therefore,
		$$( \lambda +1) x_u-2d(u) x_u=(\lambda +1) x_v-2d(v) x_v .$$
		Since $u$ and $v$ are contained in the exact same hyperedges, $d(u)=d(v)$. Hence $x_u=x_v$ if $ \lambda+1 \neq 2d(u) .$
	\end{proof}
	\begin{theorem}\label{ssofs}
 Let	$S_n$ be a star on n vertices.  Then Seidel spectrum of $S_{n}^k$ $(k\geq 2)$ is
		$$\displaystyle{\sigma_s(S_{n}^k)=\begin{pmatrix}
					1 & 3-2k & r_1 & r_2\\
					(n-1)(k-2) & n-2 & 1 & 1
				\end{pmatrix}}, $$ 
    where $r_1$ and $r_2$ are the roots of the equation, 
    $$\lambda^2-((k-1)(n-3)+1)\lambda-(n -1)(k -1) = 0.$$
   	\end{theorem} 
	\begin{proof}
		Let $e\in E(S_n)$. Adding $k-2$ vertices to $e$ forms a hyperedge of $S_n^k$. Therefore, each edge $e$ forms a hyperedge $e^k$ of $k$ vertices in $S_n^k$. Let $\{u_1,u_2,u_3,...,u_{k-1} \} \in e^k$ be vertices of degree 1. For $2\leq i \leq k-1 $ we construct $(k-2)$ linearly independent eigenvectors $\mathbf{x}^i=(x^i)_v,v\in V(S_n^k) $ as follows
		\begin{equation*}
			\mathbf{x}^i=(x^i)_v= \begin{cases}
        -1 & \text{if } v =u_1\\
        1 & \text{if } v=u_i\\
        0 & \text{otherwise.} 
    \end{cases}
		\end{equation*}
		Repeating the construction for other hyperedges, we get $(n-1)(k-2)$ linearly independent vectors associated with an eigenvalue 1.\\
		
  Let \{ $e_1,e_2,e_3,...,e_{n-1}$\} be the edges of $S_n$ and \{ $e_1^k, e_2^k, e_3^k,...,e_{n-1}^k$ \} be the hyperedges of $S_n^k$. 
		For $2\leq j\leq (n-2)$
		\begin{equation*}
			\mathbf{z}^j=(z^j)_v=\begin{cases}
                     -1 & \text{if }  v\in e_1^k \text{ and } d(v)=1 \\
                     2 & \text{if }  v \in e_j^k \text{ and } d(v)=1 \\
                     -1 & \text{if }  v \in e_{j+1}^k \text{ and } d(v)=1 \\
                      0 & \text{otherwise.} 
			\end{cases}
		\end{equation*}
		and
		\begin{equation*}
		\mathbf{z}^{n-1}=(z^{n-1})_v=\begin{cases}
                     -1 & \text{if }  v\in e_1^k \text{ and } d(v)=1 \\
                     2 & \text{if }  v \in e_{n-1}^k \text{ and } d(v)=1 \\
                     -1 & \text{if }  v \in e_{2}^k \text{ and } d(v)=1 \\
                      0 & \text{otherwise.} 
			\end{cases}
		\end{equation*}
		This construction will give $n-2$ linearly independent eigenvectors $\mathbf{z}^j$ corresponding to the eigenvalue $3-2k$.
		
		Let $(\lambda,\mathbf{x})$ be an eigenpair of $S_{n}^{k}$, we have
		$$S\mathbf{x}=\lambda\mathbf{x}.$$
	Let	$E(S_{n}) = \{e_{1},e_{2},..., e_{n-1}\}$ and $\{u_{1}^{j},u_{2}^{j},..., u_{k-1}^{j}\}$ be vertices of degree 1 in $e_{j}^{k}$ and $v$ be the vertex of degree $n-1$. Also from Lemma \ref{xu_equals_xv}, since $u_{i}^{j}$ and $u_{1}^{j}$ are contained in exactly same hyperedges  $x_{u_{i}}^{j}=x_{u_{1}}^{j}$ where $1\leq j \leq n-1$, $1\leq i \leq k-1$ and $n\geq 3$.\\
		
  By expanding $S\mathbf{x}=\lambda\mathbf{x}$, we obtain the following system of equations 
		\begin{align}
			\lambda x_{v} & = -[(k-1)x_{u_{1}^{1}}+(k-1)x_{u_{1}^{2}} + ... + (k-1)x_{u_{1}^{n-1}}], \label{eqn3.1}\\
			\lambda x_{u_{1}^{1}} & = -x_{v} - (k-2)x_{u_{1}^{1}} + (k-1)x_{u_{1}^{2}}+(k-1)x_{u_{1}^{3}} + ... + (k-1)x_{u_{1}^{n-1}}, \label{eqn3.2}\\
			\lambda x_{u_{1}^{2}} & = -x_{v} +(k-1)x_{u_{1}^{1}} - (k-2)x_{u_{1}^{2}}+(k-1)x_{u_{1}^{3}} + ... + (k-1)x_{u_{1}^{n-1}}, \label{eqn3.3}\\
			\lambda x_{u_{1}^{3}} & = -x_{v} + (k-1)x_{u_{1}^{1}} + (k-1)x_{u_{1}^{2}}-(k-2)x_{u_{1}^{3}} + ... + (k-1)x_{u_{1}^{n-1}}, \label{iam4}\\
			&\vdots \notag \\ 
			\lambda x_{u_{1}^{n-1}} & = -x_{v} + (k-1)x_{u_{1}^{1}} + (k-1)x_{u_{1}^{2}}+(k-1)x_{u_{1}^{3}} + ... - (k-2)x_{u_{1}^{n-1}}\label{3.5}.
		\end{align}
		From (\ref{eqn3.2}) and (\ref{eqn3.3}) we have,
		\begin{align*}
			x_{v} & = -\lambda x_{u_{1}^{1}} -(k-2)x_{u_{1}^{1}} + (k-1)x_{u_{1}^{2}}+(k-1)x_{u_{1}^{3}} + ... + (k-1)x_{u_{1}^{n-1}},\\
			x_{v} & = (k-1) x_{u_{1}^{1}}-\lambda x_{u_{1}^{2}} -(k-2)x_{u_{1}^{2}} + (k-1)x_{u_{1}^{3}} + ... + (k-1)x_{u_{1}^{n-1}} .
		\end{align*}
		Then, 
		\begin{equation*}
			\begin{split}
				-(\lambda + k - 2)x_{u_{1}^{1}} + & (k-1)x_{u_{1}^{2}} + (k-1)x_{u_{1}^{3}}+...+ (k-1)x_{u_{1}^{n-1}} \\
				&\quad= (k-1)x_{u_{1}^{1}}-(\lambda + k - 2)x_{u_{1}^{2}}+(k-1)x_{u_{1}^{3}} +...+ (k-1)x_{u_{1}^{n-1}}.
			\end{split}
		\end{equation*}
		On simplification, we obtain, 
		$$ (\lambda + 2k - 3)x_{u_{1}^{1}}=(\lambda + 2k - 3)x_{u_{1}^{2}}.$$
		Similarly,
		$$ (\lambda + 2k - 3)x_{u_{1}^{2}}=(\lambda + 2k - 3)x_{u_{1}^{3}}.$$
		In general, 
		$$ (\lambda + 2k - 3)x_{u_{1}^{i}}=(\lambda + 2k - 3)x_{u_{1}^{i+1}},~~i=1,2,3,\cdots, n-2.$$
		Suppose $\lambda\neq -2k+3$ then,
		$$x_{u_{1}^{1}}=x_{u_{1}^{2}}=\cdots= x_{u_{1}^{n-1}}.$$
           From (\ref{eqn3.1})-(\ref{3.5}) we obtain,
	%	So we can rewrite the system (\ref{eqn3.1})-(\ref{3.5}) into two equations as follows:
		\begin{align*}
			\begin{split}
				\lambda x_{v} & = -(k-1)(n-1)x_{u_{1}^{1}}, 
			\end{split}\\  
			\begin{split}
				\lambda x_{u_{1}^{1}} & = -x_{v} - (k-2)x_{u_{1}^{1}} + (k-1)(n-2)x_{u_{1}^{1}}.
			\end{split}
		\end{align*}
		Therefore,
		\begin{equation*}
			x_{v}  = (-\lambda-(k-2)+(k-1)(n-2))x_{u_{1}^{1}}.
		\end{equation*}
  Thus,
		%Substituting in (\ref{reds5}),
		$$(\lambda^{2}+((k-2)+(k-1)(n-2))\lambda)x_{u_{1}^{1}}=-(k-1)(n-1)x_{u_{1}^{1}}.$$
		Since $x_{u_{1}^{1}}\neq0$ we get,
		$$ \lambda^{2}-((k-1)(n-3)+1)\lambda - (k-1)(n-1) = 0 .$$
		Therefore, roots $r_{1}$, $r_{2}$ of the above equation are also eigenvalues of the hyperstar. %Since the Seidal matrix is diagonalisable, the number of eigenvalues will be equal to the number of vertex of $S_{n}^{k}$. We have,
	%	$$|V(S_{n}^{k})|=(n-1)(k-1) + 1. $$
	Thus we have all $(n-1)(k-1)+1$ eigenvalues.
	\end{proof}
	% Several researchers studied the energy of the adjacency matrix of a hypergraph. Then, the energy of the Seidel matrix is denoted by $SE(G^*)$.
 Next, we determine the Seidel energy of the hyperstar $S_n^k$. Also, we obtain a relation between the Seidel energy of $G^*$ and $G^*-v$ where $v\in V$.
	\begin{theorem}
 The Seidel energy $SE(S_n^k)$ of $S_n^k$ is,
		\begin{equation*}
			SE(S_n^k)=(n-1)(3k-5)-(2k-3)+\sqrt{(k-1)^2(n-3)^2+2(k-1)(3n-5)+1}.
		\end{equation*}
	\end{theorem}
	\begin{proof}
		The Seidel energy of $S_{n}^{k}$ is,  $SE(S_{n}^{k})=\sum_{i=1}^{(n-1)(k-1)+1} |\lambda_{i}|.$\\
  From Theorem \ref{ssofs}
		\begin{equation*}
			SE(S_n^k)=|1|(n-1)(k-2)+|3-2k|(n-2)+|r_1|+|r_2|,
		\end{equation*}
		where $r_1$ and $r_2$ are the roots of the equation  
		$ \lambda^{2}-((k-1)(n-3)+1)\lambda - (k-1)(n-1) = 0. $\\
		We can notice that
		\begin{align*}
			r_1&=\frac{(k-1)(n-3)+1+\sqrt{((k-1)(n-3)+1)^2+4(k-1)(n-1)}}{2}\geq0,\\
			r_2&=\frac{(k-1)(n-3)+1-\sqrt{((k-1)(n-3)+1)^2+4(k-1)(n-1)}}{2}\leq0.
		\end{align*}
  Therefore, $$SE(S_n^k) =(n-1)(3k-5)-(2k-3)+\sqrt{(k-1)^2(n-3)^2+2(k-1)(3n-5)+1}.$$
	\end{proof}
 
	%********************************************************************************************
	
	%%%%%%%%%%%%%%%%%%%%%%%%%%%%%%%%%%%%%%%%%%%%%%%%%%%%%%%%%%%%%%%%%%%%%%%%%%%%%%%%%%%%%%%%%%%%%%%%%%%%%%%%%%%%%%%%%%%%%%%%%%%%%%%%%%%
	For convenience in the next theorem, we denote the Seidel matrix of $G^*$ is denoted by $S(G^*)$.
	\begin{theorem}
 
		Let $G^*=(V,E)$ be a hypergraph of order $n$ and $v\in V(G^*)$ be any arbitrary vertex. Then 
		$$ SE(G^*) \geq SE(G^*-v).$$
	\end{theorem}
	\begin{proof}
		Let $\mu_1' \geq \mu_2' \geq \mu_3' \geq \cdots \geq \mu_t',$ where $t \leq n-1$ be the positive eigenvalues of $S(G^*-v)$.\\
		Then 
		$$ SE(G^*-v) = 2 \sum_{i=1}^t \mu_i'. $$
		The Seidel matrix $S(G^*-v)$ is a principal submatrix of $S(G^*)$ of order $n-1$. By Theorem \ref{interlacingthm}  we can find eigenvalues of $S(G^*)$, $\mu_1 \geq \mu_2 \geq \mu_3 \geq \cdots\geq\mu_t$ such that 
		$$ \mu_1 \geq \mu_1'\geq \mu_2 \geq \mu_2'\geq \mu_3 \geq \cdots\geq \mu_t \geq \mu_t' .$$
  Therefore,
		$$ SE(G^*)\geq 2 \sum_{i=1}^t \mu_i \geq 2 \sum_{i=1}^t \mu_i' = SE(G^*-v).$$
		Hence the result.
	\end{proof}
	\begin{corollary}
Let $S_{n}^k$ be a $k$-uniform hyperstar, then
		$$ SE(S_n^k) \geq SE(S_n^{k-1}).$$
	
	\end{corollary}
 \begin{comment}
     \begin{proof} 
		Let $S_n$ be a star graph on $n$ vertices $v_0, v_1,v_2,\cdots,v_{n-1}$ with centre vertex $v_0$ and $S_n^k$ be the $k$-uniform hypergraph on $(n-1)(k-1)+1$ vertices. By deleting the vertices $v_0, v_1,v_2,\cdots,v_{n-1}$ of $S_n^k$, we obtain $(k-1)$-uniform hyperstar $S_{n}^{k-1}$. 
		Then,
		$$ SE(S_n^k) \geq SE(S_n^{k}-v_1) \geq SE(\{S_n^{k}-v_1\}-v_2)... \geq SE(S_n^{k}-v_1-v_2-v_3\cdots-v_{n-1})=SE(S_n^{k-1}).$$
  Hence the result.
	\end{proof}
 \end{comment}
	
	\subsection{Main Seidel eigenvalues of hyperstar}
	 The main eigenvalues of a graph have been studied in  \textnormal{\cite{Cvetkovic1970}}. In this section, we discuss the number of main Seidel eigenvalues and the largest eigenvalue of a hyperstar.
  \begin{lemma}
      Let $S$ be the Seidel matrix of the hypergraph $G^*$ of order $n$. Then the rank of the matrix $\begin{bmatrix}
          \mathbfit{j} & S\mathbfit{j} & S^2\mathbfit{j}&\cdots&S^{n-1}\mathbfit{j}
      \end{bmatrix}$ is equal to the number of main Seidel eigenvalues of $G^*$.
  \end{lemma}
  \begin{proof}
      By applying similar arguments as in the proof of Theorem \ref{Hagos2002}, we get the desired result. 
  \end{proof}
	\begin{theorem}\label{thm4.9}
 The main Seidel eigenvalues of a hyperstar $S_n^k$ are $r_{1}$ and $r_{2}$ which are the roots of the equation $\lambda^2-((k-1)(n-3)+1)\lambda-(n -1)(k -1) = 0$. 
	\end{theorem}
	\begin{proof}
 \begin{comment}
     Based on Theorem \ref{ssofs} we can observe that for a hyperstar $ S_{n}^{k} $, the eigenvectors are of the form,
\begin{equation*}
\mathbf{x}^i =\left\{
   \begin{array}{ll}
     (x^i)_{u_1}=-1    &  \\
     (x^i)_{u_i}=1     &  \\
      (x^i)_{u_j}=0    & \mbox{, } u_j \in V(S_n^k-\{u_1,u_i\}) 
   \end{array}.
   \right.
\end{equation*}    
where $ 2 \leq i \leq k-1 $ and  $ \left\lbrace u_{1}, u_{2},...,u_{k-1}\right\rbrace  $ are the vertices of degree 1 of an edge $e$ . By applying the same procedure to the remaining hyperedges of $ S_{n}^{k} $, we obtain $(n-1)(k-2)$ linearly independent eigenvectors corresponding to the Seidel eigenvalue 1. Therefore, $ j^Tx^{i}=0 $ for all $(n-1)(k-2)$ eigenvectors.\\
Let $ \left\lbrace e_{1}^{k}, e_{2}^{k},...,e_{n-1}^{k}\right\rbrace  $ be edges of $ S_{n}^{k} $. Then we can construct $(n-2)$ eigenvectors corresponding to eigenvalue $3-2k$ as follows.\\ For $ 2 \leq j \leq n-2 ,$
  \begin{equation*}
			\mathbf{z}^j=\begin{cases}
				(z^j)_{u}=-1 &\text{,\;$u\in (e_1)^k$ and $d(u)=1$}\\
				(z^j)_{v}=2 & \text{,\;$v \in (e_j)^k$ and $d(v)=1$}\\
				(z^j)_{w}=-1 & \text{,\;$w \in  (e_{j+1})^k$ and $d(w)=1$}\\
				0 &\text{,\; otherwise}
			\end{cases}
		\end{equation*}
		and
		\begin{equation*}
		\mathbf{z}^{n-1}=\begin{cases}
				(z^{n-1})_{u}=-1 &\text{,\;$ u\in (e_1)^k$ and $d(u)=1$}\\
				(z^{n-1})_{v}=2 & \text{,\;$v \in (e_{n-1})^k$ and $d(v)=1$}\\
				(z^{n-1})_{w}=-1 & \text{,\;$w \in (e_2)^k$ and $d(w)=1$}\\
				0 &\text{,\; otherwise}
			\end{cases}.
		\end{equation*}
 \end{comment}
 From the proof of Theorem \ref{ssofs} we obtain $\mathbfit{j}^T\mathbf{x}^{i}=0 $ for all $(n-1)(k-2)$ eigenvectors and $ \mathbfit{j}^{T}\mathbf{z}^{j}=0 $ for all $(n-2)$ eigenvector corresponding to eigenvalue $1$ and $3-2k$ respectively.
%	Clearly $ j^{T}z^{j}=0 $ for all such $(n-2)$ eigenvector. Therefore,  1 and $3-2k$  are not main Seidel eigenvalues. 
 Thus, the only possible main Seidel eigenvalues are $ r_{1} $and $ r_{2} $, where
		$$ r_{1}=\frac{\left( k-1\right) \left( n-3\right) +1+ \sqrt{\left( \left( k-1\right) \left( n-3 \right)+1  \right)^{2}+4\left( k-1 \right) \left( n-1\right)  }}{2} $$
  and
		$$ r_{2}=\frac{\left( k-1\right) \left( n-3\right) +1- \sqrt{\left( \left( k-1\right) \left( n-3 \right)+1  \right)^{2}+4\left( k-1 \right) \left( n-1\right)  }}{2}. $$
		%Rank(W)=Rank($ W^{T} $)\\
		%Therefore suppose if $ s^{2}j=c_{1}j+c_{2}sj $\\
		%$ s^{3}j=s\left( s^{2}j\right)=s\left( c_{1}j+c_{2}sj \right)   $=$ c_{1}sj+c_{2}s^{2}j $\\
		%$ \vdots $\\
		%$ s^{n-1}j$=$ c_{1}s^{n-3}j+c_{2}s^{n-2}j $\\
		
  Next we find the rank$(\begin{bmatrix}
          \mathbfit{j} & S\mathbfit{j} & S^2\mathbfit{j}&\cdots&S^{n-1}\mathbfit{j}
      \end{bmatrix})$. Now we prove that $\mathbfit{j}$ and $S\mathbfit{j}$ are linearly independent. We can represent the Seidel matrix of $S_{n}^{k} $ as follows,
  
            \begin{equation*}
		   S=
                \begin{bmatrix}
				O_{1\times1} & -J_{1 \times \left( k-1\right) } \otimes J_{1 \times \left( n-1\right) } \\
				-J_{\left( k-1\right)\times 1 } \otimes J_{\left( n-1\right)\times 1 } & I_{n-1}\otimes B_{k-1} + (J_{n-1}-I_{n-1})\otimes J_{k-1}
			\end{bmatrix}
		\end{equation*}
		where $B_n=I_n-J_n.$
            \begin{align*}
               S\mathbfit{j}&=
			\begin{bmatrix}
				O_{1\times1} & -J_{1 \times \left( k-1\right) \left( n-1\right) } \\
				-J_{\left( k-1\right) \left( n-1\right) \times 1 } & I_{n-1}\otimes B_{k-1} + (J_{n-1}-I_{n-1})\otimes J_{k-1}
			\end{bmatrix}
			\begin{bmatrix}
				I_{1} \\
				J_{\left( n-1\right)\left( k-1\right)\times 1 }
			\end{bmatrix}\\
                &=\begin{bmatrix}
				-J_{1 \times \left( k-1\right) \left( n-1\right)} J_{\left( n-1\right)\left( k-1\right)\times1 }  \\
				-J_{\left( k-1\right)\left( n-1\right)\times 1 } +\left( I_{n-1}\otimes B_{k-1} + (J_{n-1}-I_{n-1})\otimes J_{k-1}\right) \left( J_{\left( n-1\right)\left( k-1\right)\times1 }\right)
			\end{bmatrix}.
            \end{align*}
		%\begin{equation*}
            
		%\end{equation*}
		%\begin{equation*}
             % \;\;=
			
		%\end{equation*}
		\noindent Therefore,  
		\begin{equation*}
             S\mathbfit{j}=
			\begin{bmatrix}
				-\left( n-1\right) \left( k-1\right) J_{1} \\
				-\left( n-3\right) \left( k-1\right) J_{\left( n-1\right) \left(k-1\right)\times1  }
			\end{bmatrix}.
		\end{equation*}
	Thus, $\mathbfit{j}$ and $ S\mathbfit{j} $ are linearly independent. Hence rank$(\begin{bmatrix}
          \mathbfit{j} & S\mathbfit{j} & S^2\mathbfit{j}&\cdots&S^{n-1}\mathbfit{j}
      \end{bmatrix})= 2 $. Therefore, $ r_{1} $ and $ r_{2} $ are the main Seidel eigenvalues of a hyperstar.
	\end{proof}
\begin{remark}

From Theorem \ref{thm4.9}, we can say that largest Seidel eigenvalue of hyperstar is a main Seidel eigenvalue. But the converse need not be true.   
% $r_{1}$ and $r_{2}$ are the main Seidel eigenvalues of a hyperstar $S_n^k$ Furthermore,  $r_1$ is the largest eigenvalue. Since eigenvalues $r_2$ and $2-3k$ are less than zero and $r_1$ is greater than 1. Therefore, the largest Seidel eigenvalue is a main Seidel eigenvalue.
\end{remark}	
	%%%%%%%%%%%%%%%%%%%%%%%%%%%%%%%%%%%%%%%%%%%%%%%%%%%%%%%
 \section{Spectrum of uniform double hyperstar}
 In this section, we estimate the adjacency spectrum and Seidel spectrum of uniform double hyperstar. 
 
       \begin{theorem}
               Let $S_{n_1,n_2}$ be a double star of order $n_1+n_2$, then the spectrum of $S_{n_1,n_2}^k,~(k\geq 3)$ is given by,
               $$\sigma_A( S_{n_1,n_2}^k)=\begin{pmatrix}
                 -1 & k-2 &r_1 &r_2&r_3&r_4&r_5\\
                 (k-2)(n_1+n_2-1)-1 & n_1+n_2-4&1&1&1&1&1
               \end{pmatrix}$$
               where $r_i, i=1,2,3,4,5$ are the roots of the equation
               \begin{multline}\label{eqn1}
                   \lambda^5-(-7+3 k)\lambda^4 +(17+3 k^2+n_2+n_1-k (14+n_2+n_1)) \lambda^3-(-5(3+n_2+n_1)+k(17+7n_2+7n_1)\\-k^2(7+2n_2+2n_1)+k^3)\lambda^2-(-1+(-7+12k-6k^2+k^3)n_1+(7-5k+k^2+(1-k)n_1)(-1+k)n_2)\lambda\\-(-1+k)(-1+n_2+n_1+(3-4k+k^2)n_2 n_1)=0.
               \end{multline}  
          \end{theorem}
         \begin{proof}
         Let $A(S_{n_1}^k)$ and $A(S_{n_2}^k)$ be the adjacency matrix corresponding to  $S_{n_1}^k$ and $S_{n_1}^k$ respectively. Let $D$ be a $((n_1-1)(k-1)+1) \times((n_2-1)(k-1)+1)$ matrix with the first entry equal to 1 and all other entries being 0 % $D_{((n_1-1)(k-1)+1),((n_2-1)(k-1)+1)}=\left\{ \begin{array}{rcl}
%d_{1,1}=1& \mbox{} & \\ 0 & \mbox{,} & \text{otherwise} 
% \end{array}\right.$ 
% $C_1=\left\{ \begin{array}{rcl} c_{1_{1,i}}=1& \mbox{,} & i=1,2,\cdots,k-2\\ 0 & \mbox{,} & \text{otherwise} 
% \end{array}\right.$ and 
 %$C_2=\left\{ \begin{array}{rcl}c_{2_{1,i}}=1& \mbox{,} & i=1,2,\cdots,k-2\\ 0 & \mbox{,} & \text{otherwise} 
% \end{array}\right.$
and, $C_1, C_2$ are matrices of order $\displaystyle ((n_1-1)(k-1)+1)\times(k-2), $ $((n_2-1)(k-1)+1)\times(k-2)$ respectively, with the first-row entries equal to 1 and all other entries equal to 0.
\begin{equation*}
    A(S_{n_1,n_2}^k)=\displaystyle \begin{bmatrix}
                 A(S_{n_1}^k)& D &C_1 \\
                 D^T &A(S_{n_2}^k) &C_2\\
                 C_1^T & C_2^T & J_{k-2}-I_{k-2}
                  \end{bmatrix}.
\end{equation*}
Then the characteristic polynomial of $A(S_{n_1,n_2}^k)$ is given by,
\begin{equation*}
     det(A(S_{n_{1},n_{2}}^k)-\lambda I)=det
    \begin{pmatrix}
    A(S_{n_1}^k)-\lambda I& D &C_1 \\
                 D^T & A(S_{n_2}^k)-\lambda I &C_2\\
                 C_1^T & C_2^T & J_{k-2}-(1+\lambda)I_{k-2}
   \end{pmatrix}.
   \end{equation*}
   By Lemma \ref{detblock}
   \begin{equation}\label{chareqn}
     det(A(S_{n_{1},n_{2}}^k)-\lambda I)=det(P)det(J_{k-2}-(1+\lambda)I_{k-2}) , 
   \end{equation}
   where 
   \begin{equation}\label{Pmatrix}
        P= \begin{bmatrix}
                 A(S_{n_1}^k)-\lambda I& D  \\
                 D^T &A(S_{n_2}^k)-\lambda I 
                  \end{bmatrix}-\begin{bmatrix}
                      C_1\\ C_2
                  \end{bmatrix}(J_{k-2}-(1+\lambda )I_{k-2})^{-1}\begin{bmatrix}
                      C_1^T & C_2^T
                  \end{bmatrix}.
   \end{equation}
  
   \noindent By Lemma \ref{detij} we get,$$(J_{k-2}-(1+\lambda )I_{k-2})^{-1}=\displaystyle\frac{-I_{k-2}}{1+\lambda}+\displaystyle\frac{J_{k-2}}{(1+\lambda)(k-3-\lambda)}.$$
   On simplification, we obtain
  \begin{equation}\label{P_}
      \begin{bmatrix}
                      C_1\\ C_2
                  \end{bmatrix}(J_{k-2}-(1+\lambda )I_{k-2})^{-1}\begin{bmatrix}
                      C_1^T & C_2^T
                  \end{bmatrix}=\begin{bmatrix}
                      P_1 & P_2\\
                      P_3 & P_4
                  \end{bmatrix},
  \end{equation}
  
                  where $P_1,P_2,P_3$ and $P_4$ are matrices of order $((n_1-1)(k-1)+1) \times ((n_1-1)(k-1)+1),$ $((n_1-1)(k-1)+1) \times ((n_2-1)(k-1)+1),$ $((n_2-1)(k-1)+1) \times ((n_1-1)(k-1)+1)$ and $ ((n_2-1)(k-1)+1) \times ((n_2-1)(k-1)+1)$ respectively with first entry of the matrix equal to $p$ and all other entries being zero, where $p=\displaystyle\frac{k-2}{k-3-\lambda}$ sum of all entries of $(J_{k-2}-(1+\lambda )I_{k-2})^{-1}.$ \\
                  From (\ref{Pmatrix}) and (\ref{P_}),we obtain
                   \begin{equation}\label{P_matrix}
        P= \begin{bmatrix}
                 A(S_{n_1}^k)-\lambda I-P_1& D-P_2  \\
                 D^T-P_3 &A(S_{n_2}^k)-\lambda I-P_4
                  \end{bmatrix}.
   \end{equation}
                  Let $\overline{A(S_{n_1}^k)-\lambda I}$ and $\overline{A(S_{n_2}^k)-\lambda I}$ be the matrices obtained after deleting the first row and first column of $A(S_{n_1}^k)-\lambda I $ and $A(S_{n_2}^k)-\lambda I$ respectively. Then, 
                  \begin{multline}\label{eqn5}
                      det(P)=det\left(A(S_{n_1}^k)-\lambda I-P_1\right)det \left(A(S_{n_2}^k)-\lambda I-P_4\right)\\-(1-p)^2 det\left(\overline{A(S_{n_1}^k)-\lambda I}\right)det\left(\overline{A(S_{n_2}^k)-\lambda I}\right).
                  \end{multline}
        Then, we obtain as follows
        \begin{equation}\label{eqn6}
         det\left(A(S_{n_1}^k)-\lambda I-P_1\right)=det\left(A(S_{n_1}^k)-\lambda I\right)-p \hspace{.1cm}det\left(\overline{A(S_{n_1}^k)-\lambda I}\right)   
        \end{equation}
        and
        \begin{equation}\label{eqn7}
          det\left(A(S_{n_2}^k)-\lambda I-P_4\right)=det\left(A(S_{n_2}^k)-\lambda I\right)-p \hspace{.1cm}det\left(\overline{A(S_{n_2}^k)-\lambda I}\right).  
        \end{equation}
        Also,
        \begin{align}\label{eqn8}
          det\left(\overline{A(S_{n_1}^k)-\lambda I}\right)&=det\left(I_{n_1-1} \otimes \left(J_{k-1}-(1+\lambda)I_{k-1}\right)\right)\notag\\
          &=(-\lambda-1)^{(n_1-1)(k-2)}(-\lambda+(k-2))^{n_1-1}.
        \end{align}
        Similarly, we have
        \begin{equation}\label{eqn9}
          det\left(\overline{A(S_{n_2}^k)-\lambda I}\right)=(-\lambda-1)^{(n_2-1)(k-2)}(-\lambda+(k-2))^{n_2-1}.  
        \end{equation}
        
        \noindent From Theorem \ref{Cardoso2022}, we have
        \begin{equation}\label{eqn10}
            det\left(A(S_{n_i}^k)-\lambda I\right)=(-\lambda-1)^{(n_i-1)(k-2)}(-\lambda+(k-2))^{n_i-2}(\lambda^2-(k-2)\lambda-(n_i-1)(k-1)),~i=1,2.
        \end{equation}
        From (\ref{eqn6}),(\ref{eqn8}) and (\ref{eqn10}), we obtain
        \begin{multline}\label{eqn11}
          det\left(A(S_{n_1}^k)-\lambda I-P_1\right)=\Bigl((\lambda^2-(k-2)\lambda-(n_1-1)(k-1))-\left(\frac{k-2}{k-3-\lambda}\right)(-\lambda+(k-2))\Bigr) \\ (-\lambda-1)^{(n_1-1)(k-2)}(-\lambda+(k-2))^{n_1-2}.
        \end{multline}
         %Similar equation holds for  $ det\left(A(S_{n_2}^k)-\lambda I-P_4\right)$. 
         Similarly,
     \begin{multline}\label{eqn12}
          det\left(A(S_{n_2}^k)-\lambda I-P_4\right)=\Bigl((\lambda^2-(k-2)\lambda-(n_2-1)(k-1))-\left(\frac{k-2}{k-3-\lambda}\right)(-\lambda+(k-2))\Bigr) \\ (-\lambda-1)^{(n_2-1)(k-2)}(-\lambda+(k-2))^{n_2-2}.
        \end{multline}
        Therefore,
    %Since $p=\displaystyle\frac{k-2}{k-3-\lambda}$ and substituting (\ref{eqn8}), (\ref{eqn9}) ,(\ref{eqn11}) and (\ref{eqn12}) in (\ref{eqn5}) , we obtain    
        \begin{multline}\label{eqn13}
                      det(P)=\biggl(\Bigl((\lambda^2-(k-2)\lambda-(n_1-1)(k-1))-\left(\frac{k-2}{k-3-\lambda}\right)(-\lambda+(k-2))\Bigr) \\ 
                      \Bigl((\lambda^2-(k-2)\lambda-(n_2-1)(k-1))-\left(\frac{k-2}{k-3-\lambda}\right)(-\lambda+(k-2))\Bigr) \\-\left(\frac{-\lambda-1}{k-3-\lambda}\right)^2(-\lambda+(k-2))^2 \biggr)  (-\lambda-1)^{(n_1+n_2-2)(k-2)}(-\lambda+(k-2))^{n_1+n_2-4}.
                  \end{multline}
    Since $det(J_{k-2}-(1+\lambda)I_{k-2})=(-1-\lambda)^{k-3}(k-3-\lambda)$, we get 
    %substituting (\ref{eqn13}) in (\ref{chareqn}) we get
    \begin{multline*}
                      det(P)=\biggl(\Bigl((\lambda^2-(k-2)\lambda-(n_1-1)(k-1))-\left(\frac{k-2}{k-3-\lambda}\right)(-\lambda+(k-2))\Bigr) \\ 
                      \Bigl((\lambda^2-(k-2)\lambda-(n_2-1)(k-1))-\left(\frac{k-2}{k-3-\lambda}\right)(-\lambda+(k-2))\Bigr) \\-\left(\frac{-\lambda-1}{k-3-\lambda}\right)^2(-\lambda+(k-2))^2 \biggr)  (-\lambda-1)^{(k-2)(n_1+n_2-1)-1}(-\lambda+(k-2))^{n_1+n_2-4}(k-3-\lambda).
                  \end{multline*}
    After simplification we get the desired result.

  \end{proof} 
  \begin{theorem}\label{ds}
 Let	$S_{n_1,n_2}$ be a double star on $n_1+n_2$ vertices.  Then Seidel spectrum of $S_{n_1,n_2}^k$ $(k\geq 3)$ is
		$$ \sigma_s(S_{n_1,n_2}^k)=\begin{pmatrix}
					1 & -2k+3 & r_1 & r_2 & r_3 & r_4 & r_5\\
					(k-2)(n_1+n_2-1)-1 & n_1+n_2-4 & 1 & 1 & 1 & 1 & 1
				\end{pmatrix}, $$ 
  where $r_1,~r_2,~r_3,~r_4$ and $r_5$ are the eigenvalues of the quotient matrix $Q$ of $S(S_{n_1,n_2}^k).$
  % roots of the equation $\lambda^5-\bigl(k(n_1+n_2-7)-n_1-n_2+11\bigr)\lambda^4+2(3-2k)\bigl(k(n_1+n_2-4)-n_1-n_2+7\bigr)\lambda^3+\bigl(-4 k^3*(n_1+n_2-3)+8 k^2(n_1+n_2+n_1 n_2-7)-k(4n_1+4 n_2+16 n_1 n_2-94)+8 n_1 n_2-54\bigr)\lambda^2-\bigl(8 k^3 (2n_1+2n_2-2n_1n_2-1)+k^2 (32-68n_1-68n_2+64n_1 n_2)+k (-36+90n_1+90n_2-80n_1 n_2)+11-38n_1-38n_2+32n_1 n_2\bigr)\lambda+(k-1)\bigl(4k^2(n_1+n_2-1)-4k(n_1+n_2-1+2 n_1 n_2)+n_1+n_2-1+8 n_1 n_2\bigr)=0.$
   
   	\end{theorem} 
	\begin{proof}
		Let $ E(S_{n_1})=\{e_1,e_2,e_3,\cdots,\displaystyle e_{n_1 -1}\}$ and $ E(S_{n_2})=\{e_1',e_2',e_3',\cdots,\displaystyle e_{n_2 -1}'\}$ be the edge sets  of the star graph $S_{n_1}$ and $S_{n_2}$ respectively. Let $e_0$ be the edge connecting the central vertices of $S_{n_1}$ and $S_{n_2}$. Therefore, $S_{n_1,n_2}^k$ is a $k$-uniform hypergraph obtained by adding $k-2$ vertices to every edge of $S_{n_1,n_2}$. Then $E(S_{n_1,n_2}^k)=\{e_0^k,e_1^k,e_2^k,e_3^k,\cdots,\displaystyle e_{n_1 -1}^k,e_1'^k,e_2'^k,e_3'^k,\cdots,\displaystyle e_{n_2 -1}'^k\}$. Let $\{ v_1,v_2,\cdots,v_{k-1}\}\in e_1^k $ be the vertices of degree 1. For $2\leq i \leq k-1 $ we can construct $k-2$ linearly independent eigenvectors $\mathbf{x}^i $ as follows,
  \begin{equation*}
			\mathbf{x}^i=\begin{cases}
				(x^i)_{v_1}=-1 \\
				(x^i)_{v_i}=1 \\
				(x^i)_{v_j}=0 & \text{for $v_j \in V(S_{n_1,n_2}^k-\{v_1,v_i\})$}.
			\end{cases}
		\end{equation*}
  Applying this construction on the hyperedges $e_1^k,e_2^k,e_3^k,\cdots,\displaystyle e_{n_1 -1}^k$ we obtain $(n_1-1)(k-2)$ eigenvectors corresponding to the eigenvalue $1$. By using similar construction on the hyperedges $e_1'^k,e_2'^k,e_3'^k,\cdots,$ $ e_{n_2 -1}'^k$ and on the hyperedge $e_0^k$ we can find another set of $(n_2-1)(k-2)+(k-3)$ eigenvectors associated with an eigenvalue $1$. Thus we obtain total $(k-2)(n_1+n_2-1)-1$ eigenvectors associated with the eigenvalue 1.\\
  
  For $2\leq j\leq (n_1-1),$ let $\mathbf{z}^j$ be an eigenvector corresponding to $-2k+3$ such that,
		\begin{equation*}
			\mathbf{z}^j=\begin{cases}
				(z^j)_{v_i}=1 &\text{,\;$ i=1,2,\cdots,k-1$}\\
				(z^j)_{v}=-1 & \text{,\; $v \in e_j^k$ and $d(v)=1$}\\
				0 &\text{,\; otherwise},
			\end{cases}
		\end{equation*} 
  and for $2\leq j\leq (n_2-1)$ let $\mathbf{z}_*^j$ be an eigenvector corresponding to $-2k+3$ such that,
  \begin{equation*}
			\mathbf{z}_*^j=\begin{cases}
				(z_*^j)_{v_i'}=1 &\text{,\;$v_i'\in e_1'^k$ and $d(v_i')=1$}\\
				(z_*^j)_{v'}=-1 & \text{,\; $v' \in e_j'^k$ and $d(v')=1$}\\
				0 &\text{,\; otherwise}.
			\end{cases}
		\end{equation*} 
  Therefore, we obtain $n_1+n_2-4$ eigenvectors $\mathbf{z}^j$ and $\mathbf{z}_*^j$ corresponding to the eigenvalue $-2k+3$.\\
  
  Next, we have to find the remaining eigenvalues, for that we partitioned the Seidel matrix $S(S_{n_1,n_2}^k)$ is partitioned as follows,
  \begin{equation*}
    \footnotesize \begin{bmatrix}
    0& -J_{1\times (n_1-1)(k-1)} &-1 & J_{1\times (n_2-1)(k-1)} & -J_{1\times( k-2)} \\
    -J_{(n_{1}-1)(k-1)\times 1} & I_{(n_1-1)(k-1)}+J_{(n_1-1)(k-1)} &J_{(n_{1}-1)(k-1)\times 1 } &J_{(n_1-1)( k-1)\times (n_2-1)( k-1)}& J_{(n_1-1)( k-1)\times (k-2)} \\
    &-2(I_{n_1-1}\otimes J_{k-1}) &&&\\
    -1 & J_{1\times (n_1-1)(k-1)} & 0 & -J_{1\times (n_2-1)(k-1)}& -J_{1 \times (k-2)} \\
    J_{(n_{2}-1)(k-1)\times 1} & J_{(n_2-1)( k-1)\times (n_1-1)( k-1)} & -J_{(n_{2}-1)(k-1)\times 1 } & I_{(n_2-1)(k-1)}+J_{(n_2-1)(k-1)}& J_{(n_2-1)( k-1)\times (k-2)} \\
    &&&-2(I_{n_2-1}\otimes J_{k-1}) &\\
    -J_{( k-2)\times 1} & J_{( k-2)\times (n_{1}-1)( k-1)} &-J_{(k-2) \times 1} & J_{( k-2) \times (n_{2}-1)( k-1)} & I_{k-2}-J_{k-2}
                  \end{bmatrix}.
    \end{equation*}
    Then, the quotient matrix of $S(S_{n_1,n_2}^k)$  is given by,
    \begin{equation*}
    Q=
     \begin{bmatrix}
     0 & -(n_1-1)(k-1)&-1&((n_2-1)(k-1)&2-k\\
     -1 &(n_1-3)(k-1)+1 & 1&(n_2-1)(k-1)&k-2\\
     -1&(n_1-1)(k-1)&0 &-(n_2-1)(k-1)&2-k\\
      1&(n_1-1)(k-1)&-1 &(n_2-3)(k-1)+1& k-2\\
       -1&(n_1-1)(k-1)&-1 &(n_2-1)(k-1)&3-k
     \end{bmatrix}.
    \end{equation*}
     By Theorem \ref{equipartition}, the spectrum of $S(S_{n_1,n_2}^k)$ contains the spectrum of $Q$. Thus we have all $(n_1+n_2-3)k$ eigenvalues.
     
     % Therefore the characteristic polynomial of Q, $\lambda^5-\bigl(k(n_1+n_2-7)-n_1-n_2+11\bigr)\lambda^4+2(3-2k)\bigl(k(n_1+n_2-4)-n_1-n_2+7\bigr)\lambda^3+\bigl(-4 k^3*(n_1+n_2-3)+8 k^2(n_1+n_2+n_1 n_2-7)-k(4n_1+4 n_2+16 n_1 n_2-94)+8 n_1 n_2-54\bigr)\lambda^2-\bigl(8 k^3 (2n_1+2n_2-2n_1n_2-1)+k^2 (32-68n_1-68n_2+64n_1 n_2)+k (-36+90n_1+90n_2-80n_1 n_2)+11-38n_1-38n_2+32n_1 n_2\bigr)\lambda+(k-1)\bigl(4k^2(n_1+n_2-1)-4k(n_1+n_2-1+2 n_1 n_2)+n_1+n_2-1+8 n_1 n_2\bigr)$ is a factor of characteristic polynomial of $S(S_{n_1,n_2}^k)$. Hence the result.
  \end{proof}
  %%%%%%%%%%%%%%%%%%%%%%%%%%%%%%%%%%%%%%%%%%%%%%%%%%%%%%%%___________________________________________________

\section{Spectrum of sunflower hypergraph }
In this section, we estimate the adjacency and Seidel eigenvalues of sunflower hypergraph.
	\begin{theorem}
	    Let $S^k$ be a $k$-uniform sunflower hypergraph. If $k\geq 2$ is an integer, then the characteristic
polynomial of $S^k$ is
\begin{multline}
        P_{A}(\lambda)=(1+\lambda)^{(k-1)(k-2)}\left((2-3k+k^2)+(6-6k+k^2)\lambda+(4-2k)\lambda^2+\lambda^3\right)\\
        \left((-3+2k)+(k-3)\lambda-\lambda^2\right)^{k-2}.
    \end{multline}
        \end{theorem}
        \begin{proof}
           Let $\mathbf{e}_{i}\in \mathbb{R}^k$  with one in the $i$-th coordinate and zero elsewhere. Then, the adjacency matrix of $S^k$ can be written as
           $$A(S^k)=\begin{bmatrix}
              J_k-I_k  & \mathbf{e}_{2}\otimes J_{1,k-1} &  \mathbf{e}_{3}\otimes J_{1,k-1} &\cdots&  \mathbf{e}_{k}\otimes J_{1,k-1}  \\
              \mathbf{e}^{T}_{2}\otimes J_{k-1,1}  & J_{k-1}-I_{k-1} & \mathbf{0}_{k-1} & \cdots & \mathbf{0}_{k-1} \\
              \mathbf{e}^{T}_{3}\otimes J_{k-1,1}  & \mathbf{0}_{k-1} & J_{k-1}-I_{k-1} & \cdots & \mathbf{0}_{k-1}\\
               \vdots  & \vdots & \vdots & \ddots & \vdots\\
                \mathbf{e}^{T}_{k}\otimes J_{k-1,1}  & \mathbf{0}_{k-1} & \mathbf{0}_{k-1} & \cdots & J_{k-1}-I_{k-1}
                \end{bmatrix}.$$
         Therefore the characteristic polynomial of $S^k$ is given by
         $$ det(A(S^k)-\lambda I)=det\begin{pmatrix}
             J_k-(1+\lambda)I_k  & B\\
             B^{T}               & I_{k-1} \otimes(J_{k-1}-(1+\lambda)I_{k-1})
                              \end{pmatrix}$$
            where $B=\begin{bmatrix}
                \mathbf{e}_2\otimes J_{1,k-1} &  \mathbf{e}_{3}\otimes J_{1,k-1} &\cdots&  \mathbf{e}_{k}\otimes J_{1,k-1}  
            \end{bmatrix}.$ From Lemma \ref{detblock}, we obtain
     \begin{equation}\label{eqns1}
         det(A(S^k)-\lambda I)=det(J_k-(1+\lambda)I_k) det((I_{k-1} \otimes(J_{k-1}-(1+\lambda)I_{k-1}))-B^{T}(J_k-(1+\lambda)I_k)^{-1}B).
     \end{equation}
     Applying Lemma \ref{detij}, we have
     \begin{equation}
         (J_k-(1+\lambda)I_k)^{-1}=\frac{-I_k}{1+\lambda}+\frac{J_k}{(1+\lambda)
         (k-1-\lambda)}.
     \end{equation}
     Then,
     %After performing the usual matrix multiplication, we obtain
     \begin{align*}
         B^T(J_k-(1+\lambda)I_k)^{-1}B&=\left(\frac{-I_{k-1}}{1+\lambda}+\frac{J_{k-1}}{(1+\lambda)
         (k-1-\lambda)}\right)\otimes J_{k-1}\notag\\
         &=\frac{-I_{k-1}\otimes J_{k-1}}{1+\lambda}+\frac{J_{(k-1)^2}}{(1+\lambda)
         (k-1-\lambda)}.
     \end{align*}
  %  Since $I_{k-1} \otimes(J_{k-1}-(1+\lambda)I_{k-1})=I_{k-1} \otimes J_{k-1}-(1+\lambda)I_{(k-1)^2}$, using equation(\ref{3}) we obtain
  Therefore,
    \begin{align*}
    \begin{split}
        det((I_{k-1} \otimes(J_{k-1}- &(1+\lambda)I_{k-1}))-B^{T} (J_k-(1+\lambda)I_k)^{-1}B)\\
         &=det\left(I_{k-1} \otimes J_{k-1}-(1+\lambda)I_{(k-1)^2}+\frac{I_{k-1}\otimes J_{k-1}}{1+\lambda}-\frac{J_{(k-1)^2}}{(1+\lambda)(k-1-\lambda)}\right)\\
         &=det\left(-(1+\lambda)I_{(k-1)^2}+ \left(\frac{2+\lambda}{1+\lambda}\right)I_{k-1}\otimes J_{k-1}-\frac{J_{(k-1)^2}}{(1+\lambda)(k-1-\lambda)}\right).
    \end{split}
    \end{align*}
    % Now take $\beta=-(1+\lambda)$, $M=\left(\frac{-2-\lambda}{1+\lambda}\right)I_{k-1}\otimes J_{k-1}$ and $\rgamma=\displaystyle\frac{1}{(1+\lambda)(k-1-\lambda)}$,
    Using Lemma \ref{coronal}, we get
    \begin{align}\label{eqn4}
    \begin{split}
        det&((I_{k-1}  \otimes(J_{k-1}- (1+\lambda)I_{k-1}))-B^{T} (J_k-(1+\lambda)I_k)^{-1}B)\\
         &=\left(1-\frac{1}{(1+\lambda)(k-1-\lambda)}\rchi_M(-1-\lambda)\right)det\left((-1-\lambda)I_{(k-1)^2}+\frac{2+\lambda}{1+\lambda}I_{k-1}\otimes J_{k-1}\right),
    \end{split}
    \end{align} 
    where $M=\left(\frac{-2-\lambda}{1+\lambda}\right)I_{k-1}\otimes J_{k-1}$ and $\rchi_M(-1-\lambda)=\displaystyle\frac{(k-1)^2(1+\lambda)}{-(1+\lambda)^2+(k-1)(2+\lambda)}$.
    Since $\displaystyle\frac{(k-1)(2+\lambda)}{1+\lambda}$ and $0$ are the eigenvalues of $\displaystyle\frac{2+\lambda}{1+\lambda}I_{k-1}\otimes J_{k-1}$ with multiplicity $k-1$ and $(k-1)(k-2)$ respectively, we get
    $$det\left((-1-\lambda)I_{(k-1)^2}+\frac{2+\lambda}{1+\lambda}I_{k-1}\otimes J_{k-1}\right)=\left((k-1)(2+\lambda)-(1+\lambda)^2\right)^{(k-1)} (-1-\lambda)^{(k-1)(k-3)}.$$
    From (\ref{eqn4}) we obtain,
    \begin{align}\label{eqns5}
        det((I_{k-1} \otimes & (J_{k-1}-(1+\lambda)I_{k-1}))-B^{T}(J_k-(1+\lambda)I_k)^{-1}B)\notag\\
         &=\frac{
   \splitfrac{\left(\left(-(1+\lambda)^2+(k-1)(2+\lambda)\right)(k-1-\lambda)-(k-1)^2\right)(-1-\lambda)^{(k-1)(k-3)}}%
             {((k-1)(2+\lambda)-(1+\lambda)^2)^{k-2}}}%
             {k-1-\lambda}.
    \end{align} 
    Clearly $det(J_k-(1+\lambda)I_k)=(k-1-\lambda)(-1-\lambda)^{k-1}$. Therefore from (\ref{eqns5}) and (\ref{eqns1}) , we get
    \begin{multline*}
         det(A(S^k)-\lambda I)=(1+\lambda)^{(k-1)(k-2)}\left((k-1-\lambda)\left((k-1)(k-2)-(1+\lambda)^2\right)-(k-1)^2\right)\\ \left((k-1)(2+\lambda)-(1+\lambda)^2\right)^{k-2}.
    \end{multline*}
    On simplification, the characteristic polynomial of $S^k$ becomes
    \begin{multline*}
        P_{A(S^k)}(\lambda)=(1+\lambda)^{(k-1)(k-2)}\left((2-3k+k^2)+(6-6k+k^2)\lambda+(4-2k)\lambda^2+\lambda^3\right)\\
        \left((-3+2k)+(k-3)\lambda-\lambda^2\right)^{k-2}.
    \end{multline*}
    
        \end{proof}
    \begin{corollary}
        Let $S^k(k\geq 2)$ be a $k$-uniform sunflower hypergraph. Then the spectrum $\sigma_A(S^k)$ is
     $$ \sigma_A(S^k) =\begin{pmatrix}
               -1 & \displaystyle\frac{(k-3)+\sqrt{(k+3)(k-1)}}{2}  & \displaystyle\frac{(k-3)-\sqrt{(k+3)(k-1)}}{2} & r_{1} & r_{2} &r_{3}\\
               (k-1)(k-2) & k-2 & k-2 & 1  & 1  & 1
             \end{pmatrix},$$
           where $ r_i=\displaystyle\frac{2}{3}(-2+k+\displaystyle\sqrt{-2+2k+k^2)}~\cos\left(\displaystyle\frac{\theta+2(i-1)\pi}{3} \right)$ ,$~i=1,2,3$ and \\$\theta=\cos^{-1}\left(\displaystyle \frac{34-51k+21k^2-2k^3}{2\sqrt{(-2+2k+k^2)^3}}\right).$
    \end{corollary}
    \begin{proof}
        The characteristic polynomial of $S^k$ is given by,
    \begin{multline*}
        P_{A(S^k)}(\lambda)=(1+\lambda)^{(k-1)(k-2)}\left((2-3k+k^2)+(6-6k+k^2)\lambda+(4-2k)\lambda^2+\lambda^3\right)\\
        \left((-3+2k)+(k-3)\lambda-\lambda^2\right)^{k-2}
    \end{multline*}
    Clearly, $\displaystyle\frac{(k-3)\pm \sqrt{(k+3)(k-1)}}{2}$ are the roots of the equation $(-3+2k)+(k-3)\lambda-\lambda^2=0$. Using the method in \cite{cubicpoly}, we get $ r_i=\displaystyle\frac{2}{3}(-2+k+\sqrt{-2+2k+k^2)}\cos\left(\frac{\theta+2(i-1)\pi}{3} \right)$ , $i=1,2,3$  are the solution of the equation $(2-3k+k^2)+(6-6k+k^2)\lambda+(4-2k)\lambda^2+\lambda^3=0$ where $\theta=\cos^{-1}\left(\displaystyle \frac{34-51k+21k^2-2k^3}{2\sqrt{(-2+2k+k^2)^3}}\right).$
    \end{proof}
     \begin{theorem}
        Let $S^k~(k\geq 2)$ be a $k$-uniform sunflower hypergraph.  Then the Seidel spectrum $\sigma_S(S^k)$ of $S^k$ is
     $$ \sigma_S(S^k) =\begin{pmatrix}
               1 & \displaystyle 2-k+\sqrt{(k+3)(k-1)}  & \displaystyle 2-k-\sqrt{(k+3)(k-1)} & r_{1} & r_{2} &r_{3}\\
               (k-1)(k-2) & k-2 & k-2 & 1  & 1  & 1
             \end{pmatrix},$$
           where $r_1,r_2$ and $r_3$ are the roots of the equation $\lambda^3 -(6-5 k+k^2) \lambda^2-(-17+26k-12k^2+2k^3) \lambda-8+17 k-11k^2+2k^3=0. $
    \end{theorem}
    \begin{proof}
    Let $\eta_i\in \mathbb{R}^k$ with $-1$ in the $i$-th cordinate and $1$ elsewhere. Then the Seidel matrix of the sunflower hypergraph can be expressed as follows,
    $$S(S^k)=\begin{bmatrix}
              I_k-J_k  & \eta_{2}\otimes J_{1,k-1} &  \eta_{3}\otimes J_{1,k-1} &\cdots&  \eta_{k}\otimes J_{1,k-1}  \\
              \eta^{T}_{2}\otimes J_{k-1,1}  & I_{k-1}-J_{k-1} & J_{k-1} & \cdots & J_{k-1} \\
              \eta^{T}_{3}\otimes J_{k-1,1}  & J_{k-1,1} & I_{k-1}-J_{k-1} & \cdots & J_{k-1,1}\\
               \vdots  & \vdots & \vdots & \ddots & \vdots\\
                \eta^{T}_{k}\otimes J_{k-1,1}  & J_{k-1,1} & J_{k-1,1} & \cdots & I_{k-1}-J_{k-1}
                \end{bmatrix}.
                $$
        From Theorem \ref{seideleig}, the three Seidel eigenvalues of $S^k$ are $1,~2-k+\sqrt{(k+3)(k-1)}$  and $2-k-\sqrt{(k+3)(k-1)}$. Now, to find the multiplicity of these eigenvalues, we construct the corresponding linearly independent eigenvectors.\\
        \indent For eigenvalue 1, we determine $(k-1)(k-2)$ linearly independent eigenvectors $\mathbf{y}^i_j= \bigl(y_j^i\bigr)_v,~v\in V(S^k)$ as follows
       $$ \bigl(y_j^i\bigr)_v=\begin{cases}
	               1, & \text{if $v=v_{i,2}$},\\
                      -1, & \text{if $v=v_{i,j+1}$},\\
                       0, & \text{otherwise},
                       \end{cases}$$
        where $1\leq i \leq k-1, ~2\leq j \leq k-1$. Thus $1$ is an eigenvalue of $S^k$ with multiplicity $(k-1)(k-2).$ \\
        
        Let $\mathbf{x}^i=\begin{bmatrix}
            \mathbfit{x}^i_{1} &   \mathbfit{x}^i_{2} \cdots \mathbfit{x}^i_{k-1} 
       \end{bmatrix}^T, ~(2\leq i \leq k-1) $ be the eigenvectors corresponding to the eigenvalue $2-k+\sqrt{(k+3)(k-1)}$ of $S^k$, where $\mathbfit{x}^i_1=\begin{bmatrix}
          x^i_{v_{0,0}}  &  x^i_{v_{1,1}} &  x^i_{v_{2,1}} &   x^i_{v_{3,1}} \cdots x^i_{v_{k-1,1}}   
       \end{bmatrix}^T $ and 
       $\mathbfit{x}^i_{j+1}=\begin{bmatrix}
            x^i_{v_{j,2}} & x^i_{v_{j,3}} \cdots x^i_{v_{j+1,k}} 
       \end{bmatrix} $, $1\leq j \leq k-2$. From Lemma \ref{xu_equals_xv}, we get $\mathbfit{x}^i_{j+1}$'s are of the form $ c J_{k-1,1} $ where $c$ is any constant. Then,
       % \hl{So, we can construct $(k-2)$ eigenvectors $\mathbf{x}^i$ as follows, }
       $$\mathbf{x}^i=\begin{bmatrix}
            \mathbfit{x}^i_{1} & c_2J_{k-1,1} &   c_3J_{k-1,1 } &\cdots &c_{k-1}J_{k-1,1}  
       \end{bmatrix}^T. $$
       For $v\in \{v_{0,0},~v_{1,1},~v_{2,1},~v_{3,1},\ldots,v_{k-1,1}\}$,  $2\leq r \leq k-1$
       \begin{center}
            $ \mathbfit{x}^i_{1}=(x^i_{1})_{v}=\begin{cases}
	               1, & \text{if $v=v_{1,1}$},\\
                      -1, & \text{if $v=v_{i,1}$},\\
                       0, & \text{otherwise},
                \end{cases}$ 
                 and $c_{r}=\begin{cases}
	              \frac{1}{2}(1-\sqrt{\frac{k+3}{k-1}}), & \text{if $r=2$},\\
                      -\frac{1}{2}(1-\sqrt{\frac{k+3}{k-1}}), & \text{if $r=i+1$},\\
                       0, & \text{otherwise}.
                \end{cases}$
       \end{center}
       Therefore, we obtain a family of $k-2$ linearly independent eigenvectors associated with an eigenvalue  $2-k+\sqrt{(k+3)(k-1)}$. Hence $2-k+\sqrt{(k+3)(k-1)}$ is an eigenvalue of multiplicity $(k-2)$.\\
       
       Similarly, we can determine a set of linearly independent eigenvectors $\mathbf{z}^i~(2 \leq i \leq k-2)$ associated with an eigenvalue $2-k-\sqrt{(k+3)(k-1)}$ as follows,
       $$\mathbf{z}^i=\begin{bmatrix}
          \mathbfit{z}^i_{1} & c_2J_{k-1,1} &   c_3J_{k-1,1 } &\cdots &c_{k-1}J_{k-1,1}  
       \end{bmatrix} ^T.$$
       For $v\in \{v_{0,0},~v_{1,1},~v_{2,1},~v_{3,1},\ldots,v_{k-1,1}\}$, $2\leq r \leq k-1$  
       $$ \mathbfit{z}^i_{1}=(z^i_{1})_{v}=\begin{cases}
	               1, & \text{if $v=v_{1,1}$},\\
                      1, & \text{if $v=v_{i,1}$,}\\
                     -2, & \text{if $v=v_{i+1,1}$},\\
                       0, & \text{otherwise,}
                \end{cases} 
                 \text{and}~~ c_{r}=\begin{cases}
	              \frac{1}{2}(1+\sqrt{\frac{k+3}{k-1}}), & \text{if $r=2$},\\
                    \frac{1}{2}(1+\sqrt{\frac{k+3}{k-1}}), & \text{if $r=i+1$},\\
                    -(1+\sqrt{\frac{k+3}{k-1}}), & \text{if $r=i+2$},\\
                       0, & \text{otherwise}
                \end{cases}$$\\
                and\\
       $$\mathbfit{z}^{k-1}_{1}=(z^{k-1}_{1})_{v}=\begin{cases}
	               1, & \text{if  $v=v_{1,1}$,}\\
                     -2, & \text{if $v=v_{2,1}$},\\
                      1, & \text{if $v=v_{k-1,1}$},\\
                      0, & \text{otherwise},
                \end{cases}
                 \text{and}~~ c_{r}=\begin{cases}
	              \frac{1}{2}(1+\sqrt{\frac{k+3}{k-1}}), & \text{if $r=2$},\\
                    -(1+\sqrt{\frac{k+3}{k-1}}), & \text{if $r=3$},\\
                    \frac{1}{2}(1+\sqrt{\frac{k+3}{k-1}}), & \text{if $r=k-1$},\\
                    0, & \text{otherwise}.
                \end{cases}$$ 
        Since the eigenvectors are linearly independent, $2-k-\sqrt{(k+3)(k-1)}$ is an eigenvalue of multiplicity $k-2$.  
    The remaining eigenvalues of $S(S^k)$ are those of its quotient matrix $Q$ of $S(S^k)$,
$$Q=\begin{bmatrix}
    0 & 1-k &(k-1)^2\\
    -1& 2-k &(k-1)(k-3)\\
    1&  k-3 &(k-2)^2
\end{bmatrix}.$$
Thus the characteristic equation of $Q$ is given by,\\
$$\lambda^3-(6-5 k+k^2) \lambda^2-(-17+26k-12k^2+2k^3) \lambda-8+17 k-11k^2+2k^3=0$$
Hence, the theorem follows.
\end{proof}

 \section{Conclusion}   
 In this paper, we determine the relation between the characteristic polynomial of Seidel and the adjacency matrix of the hypergraph. In addition,  we obtain the Seidel spectrum and the number of walks of length $l$ of $(k,r)$-regular hypergraph. Also, we discuss the Seidel spectrum, Seidel energy and main Seidel eigenvalues of hyperstar. Using the adjacency matrix of hyperstar we determine the adjacency spectrum and Seidel spectrum of uniform double hyperstar. Moreover, we estimate the adjacency and Seidel spectrum of the sunflower hypergraph.

\section{Declarations}
On behalf of all authors, the corresponding author states that there is no conflict of interest.
\bibliographystyle{plain}
\bibliography{staref}
\end{document}